\documentclass[a4paper,12pt]{amsart}
\usepackage[]{amscd,amsfonts,amsmath,amsxtra,amssymb,bm}
\usepackage{graphics}
\usepackage{extarrows}
\usepackage[all]{xy}
\usepackage{hyperref}
\usepackage{eucal}
\usepackage{rotating}
\usepackage[T1]{fontenc}
\newtheorem{sub}{}[section]
\newtheorem{subsub}{}[sub]
\usepackage{rotating}
\usepackage[active]{srcltx}
\usepackage{bbm}
\def\ov#1{\overline{#1}}

\def\coker{\mathop{\rm coker}\nolimits}
\def\Hom{\mathop{\rm Hom}\nolimits}

\def\HHom{\mathop{\mathcal Hom}\nolimits}

\def\Ext{\mathop{\rm Ext}\nolimits}
\def\EExt{\mathop{\mathcal Ext}\nolimits}

\def\Pic{\mathop{\rm Pic}\nolimits}

\def\Aut{\mathop{\rm Aut}\nolimits}

\def\End{\mathop{\rm End}\nolimits}

\def\EEnd{\mathop{\mathcal End}\nolimits}

\def\imm{\mathop{\rm im}\nolimits}

\def\rk{\mathop{\rm rk}\nolimits}

\def\spec{\mathop{\rm spec}\nolimits}

\def\lra{\longrightarrow}

\def\sigg{\mathop{\hbox{$\displaystyle\sum$}}\limits}

\def\hfl#1#2{\smash{\mathop{\ \hbox to 12mm{\rightarrowfill}}
\limits^{\scriptstyle#1}_{\scriptstyle#2} \ }}
\def\hflb#1#2{\smash{\mathop{\hbox to 12mm{\leftarrowfill}}
\limits^{\scriptstyle#1}_{\scriptstyle#2}}}

\def\m#1{{\hbox{$#1$}}}

\def\ot{\otimes}
\def\og{\leavevmode\raise.3ex\hbox{$\scriptscriptstyle\langle\!\langle$}}
\def\fg{\leavevmode\raise.3ex\hbox{$\scriptscriptstyle\,\rangle\!\rangle$}}

\def\nsp{\lbrace 0\rbrace}

\def\dsp{\displaystyle}
\def\Ssect#1#2{\pagebreak[3]\begin{sub}\label{#2}{\sc\small\small
#1}\rm\medskip}

\def\sepsec{\vskip 1.4cm}
\def\sepsub{\vskip 0.7cm}
\def\sepsubsub{\vskip 0.5cm}
\def\sepprop{\vskip 0.5cm}

\def\xmat#1{\[\xymatrix{#1}\]}
\def\flinc{\ar@{^{(}->}}
\def\fleq{\ar@{=}}
\def\flon{\ar@{->>}}
\def\fmaps{\ar@{|-{>}}}
\def\fflat{\ar@{-}}
\def\fimpl{\ar@{=>}}
\def\Nligne{\hfil\break}

\def\ED{\vskip 1cm\end{document}}

\def\BS#1{{\boldsymbol{#1}}}

\newcommand{\M}{{\mathbb M}}
\newcommand{\A}{{\mathbb A}}
\newcommand{\B}{{\mathbb B}}
\newcommand{\Z}{{\mathbb Z}}

\newcommand{\C}{{\mathbb C}}

\renewcommand{\P}{{\mathbb P}}

\newcommand{\F}{{\mathbb F}}
\newcommand{\E}{{\mathbb E}}
\newcommand{\G}{{\mathbb G}}

\renewcommand{\L}{{\mathbb L}}

\newcommand{\kd}{{\mathcal D}}
\newcommand{\ke}{{\mathcal E}}
\newcommand{\kf}{{\mathcal F}}
\newcommand{\kg}{{\mathcal G}}
\newcommand{\kh}{{\mathcal H}}
\newcommand{\ki}{{\mathcal I}}

\newcommand{\kk}{{\mathcal K}}
\newcommand{\kl}{{\mathcal L}}

\newcommand{\kn}{{\mathcal N}}
\newcommand{\ko}{{\mathcal O}}

\newcommand{\kx}{{\mathcal X}}

\newcommand{\mycomment}[1]{}
\newcommand{\ie}{i.\,e.\ }
\begin{document}

\def\refname{References}
\def\contentsname{Summary}
\def\proofname{Proof}
\def\abstractname{Resume}

\author{Jean--Marc Dr\'{e}zet}
\address{
Institut de Math\'ematiques de Jussieu - Paris Rive Gauche\\
F-75252 Paris, France  --- now retired}
\email{jean-marc.drezet@wanadoo.fr}
\title[{Quasi locally free sheaves}] {Quasi locally free sheaves on primitive
double schemes}

\begin{abstract}
Let $\BS{\kx}$ be a primitive double scheme, $X=\BS{\kx}_{red}$, \m{L=\ki_X}
the ideal sheaf of $X$, which is a line bundle on $X$. Let $\ke$
be a coherent sheaf on $\BS{\kx}$. Then there exists a nonempty open subset
$U$ of $\BS{\kx}$ and integers \m{m_1}, \m{m_2}, uniquely determined by $\ke$,
such that $\ke_{|U}$ is locally isomorphic to \
\m{(\ko_X\ot\C^{m_1})\oplus(\ko_{\BS{\kx}}\ot\C^{m_2})}. If we can take
\m{U=\BS{\kx}}, we say the $\ke$ is {\em quasi locally free}. If \m{m_1=1} we
say that $\ke$ is {\em rigid}. In this case, we show that deformations of
$\ke$ are of the same type. So locally free sheaves and rigid quasi locally
free sheaves on $\BS{\kx}$ play the same role on $\BS{\kx}$ as locally free
sheaves do on smooth schemes.

Let $\ke$ be a quasi locally free sheaf on $\BS{\kx}$. Then \m{\ke_1=\ki_X\ke}
and \ \m{\ke_{|X}} \ are vector bundles on $X$, and we have a canonical
surjective morphism \ \m{\BS{\Psi}_\ke:\ke_{|X}\to\ke_1\ot L^*}. Similarly
there is a canonical exact sequence \m{\Sigma(\ke)} of length 4 of vector
bundles on $X$ associated to
$\ke$. Given a surjective morphism \ \m{\Phi:\E\to E} \ of vector bundles on
$X$, we define an element \ \m{\nabla_{\E,E}(\Phi)\in\Ext^2_X(\E,E\ot L)} \
such that there exists a quasi locally free sheaf $\ke$ on $\BS{\kx}$ such
that \ \m{\BS{\Psi}_\ke\simeq\Phi} \ if and only \ \m{\nabla_{\E,E}(\Phi)=0},
and prove that if $\E$ or $E$ can be extended to a vector bundle on $\BS{\kx}$,
then $\ke$ exists. Two ways are given to define \m{\nabla_{\E,E}}, a direct
one and a more complicated one, using the cohomology theory of sheaves of
non abelian groups, that can also be used to define an obstruction to the
existence of a quasi locally free sheaf $\ke$ corresponding to a given exact
sequence of length 4.

We study also the existence and properties of flat families of quasi locally
free sheaves, giving the possibility to construct universal families.
\end{abstract}

\maketitle
\tableofcontents


\newpage

\section{Introduction}\label{intro}

A {\em primitive multiple scheme} is a Cohen-Macaulay scheme $\BS{\kx}$ over
$\C$ such that:
\begin{enumerate}
\item[--] $\BS{\kx}_{red}=X$ is a smooth connected variety,
\item[--] for every closed point \m{x\in X}, there exists a neighborhood $U$
of $x$ in $\BS{\kx}$, and a smooth variety $S$ of dimension \ \m{\dim(X)+1}
such that $U$ is isomorphic to a closed subscheme of $S$.
\end{enumerate}
We call $X$ the {\em support} of $\BS{\kx}$. It may happen that $\BS{\kx}$ is
quasi-projective, and in this case it is projective if $X$ is.

Let \m{\ki_X} denote the ideal sheaf of $X$ in $\BS{\kx}$. We will suppose
here that $\BS{\kx}$ is a {\em primitive double scheme}, \ie \
\m{\ki_X\not=0} \ and \ \m{\ki_X^2=0}. In this case \ \m{L=\ki_X} is a line
bundle on $X$ that we call it {\em the line bundle on $X$ associated to
$\BS{\kx}$}. The {\em trivial double scheme} with induced smooth variety
$X$ and associated line bundle $L$ on $X$ is the second infinitesimal
neighborhood of $X$, embedded by the zero section in  the dual bundle $L^*$,
seen as a smooth variety.

The primitive multiple curves where defined in \cite{fe}, \cite{ba_fo}.
Primitive double curves were parameterized and studied in \cite{ba_ei} and
\cite{ei_gr}. More results on primitive multiple curves can be found in
\cite{dr2}, \cite{dr1}, \cite{dr4}, \cite{dr5}, \cite{dr6}, \cite{dr7},
\cite{dr8}, \cite{dr9}, \cite{ch-ka}, \cite{sa1}, \cite{sa-vi}, \cite{sa3}.
Some primitive double schemes are studied in \cite{b_m_r}, \cite{ga-go-pu} and
\cite{gonz1}. The case of varieties of any dimension is studied in
\cite{dr10}. In \cite{dr11}, the construction and properties of fine moduli
spaces of vector bundles on primitive multiple schemes are described, and in
\cite{dr12} more general sheaves are studied.

We have a canonical exact sequence on $X$
\[ 0\lra L=\ki_X\lra\Omega_{\BS{\kx}|X}\lra\Omega_X\lra 0 \ , \]
corresponding to \ \m{\mu\in\Ext^1_X(\Omega_X,L)=H^1(X,T_X\ot L)}. Let \
\[\sigma(\BS{\kx}) \ = \ \C\lambda\in\P(H^1(X,T_X\ot L))\cup\nsp \ . \]
If \m{\BS{\kx}'} is another primitive double scheme with reduced subscheme $X$
and associated line bundle $L$, there exists an isomorphism
\m{\BS{\kx}\to\BS{\kx}'} inducing the identity on $X$ if and only if \
\m{\sigma(\BS{\kx}')=\sigma(\BS{\kx})}. Moreover, for any \
\m{\alpha\in\P(H^1(X,T_X\ot L))\cup\nsp}, there exists a primitive double
scheme $\BS{\kx}$ such that \ \m{\sigma(\BS{\kx})=\alpha} (cf. \cite{ba_ei},
\cite{dr10}, \ref{PMS}).

The Zariski spaces of $X$ and $\BS{\kx}$ are the same. If \m{U\subset X} is an
open subset, let \m{U^{(2)}} denote the corresponding open subscheme of
$\BS{\kx}$.

Let $\ke$ be a coherent sheaf on $\BS{\kx}$. Then there is a nonempty open
subset \m{U\subset X} and uniquely determined integers \m{m_1,m_2}, such
that on \m{U^{(2)}}, $\ke$ is locally isomorphic to \Nligne
\m{(\ko_X\ot\C^{m_1})\oplus(\ko_{\BS{\kx}}\ot\C^{m_2})}. The sum \ \m{m_1+m_2}
\ is called the {\em generalized rank} of $\ke$. It is invariant by
deformation of $\ke$.

If we can take \m{U=X}, we say that $\ke$ is {\em quasi locally free}, and
\m{(m_1,m_2)} is called the {\em type} of $\ke$. If $\ke$ is quasi locally
free and \ \m{m_1=1} \ we say that $\ke$ is of {\em rigid type}.

The subject of this paper is the study of quasi locally free sheaves, their
construction from basic data on $X$, their parameterization, and the
properties of flat families of such sheaves.

\sepprop

{\bf Terminology and notations: } -- A {\em scheme} is a noetherian separated
scheme over $\C$.

-- If $Z$ is a scheme and \ \m{Y\subset Z} \ is a closed subscheme,
\m{\ki_{Y,Z}} or \m{\ki_Y} denotes the ideal sheaf of $Y$.

-- Let $Y$, $Z$ be schemes, $\ke$ a coherent sheaf on \m{Y\times Z}, and
\m{z\in Z} a closed point. Then \m{\ke_z} denotes the sheaf
\m{\ke_{|{Y\times\{z\}}}} on $Y$.

-- Let $Y$, $Z$, $T$ be schemes, $\ke$ a coherent sheaf on \m{Y\times T}, and \
\m{f:Z\to T} \ a morphism. Then \ \m{f^\sharp(\ke)=(I_Y\times f)^*(\ke)}.

-- If $M$, $N$ are finite dimensional complex vector spaces, \m{\kl(M,N)}
denotes the space of linear maps \ \m{M\to N}.

\sepsub

\Ssect{Construction of quasi locally free sheaves}{Intro1}

\begin{subsub} Associated sheaves on $X$ -- \rm Let $\ke$ be a coherent sheaf
on \m{\BS{\kx}}. Then we have the following sheaves on $X$:
\begin{enumerate}
\item[--] $\ke_{|X}$, the restriction of $\ke$ to $X$.
\item[--] $\ke_1=\ki_X\ke$ . We have \ \m{\ke_{|X}=\ke/\ke_1} .
\item[--] $\ke^{(1)}$, the subsheaf of $\ke$ of elements annihilated by
\m{\ki_X}.
\end{enumerate}
From \cite{dr10}, $\ke$ is quasi locally free if and only \m{\ke_{|X}} and
\m{\ke_1} are locally free on $X$.

The sheaf \m{\ke^{(1)}} contains \m{\ke_1}. The canonical morphism \
\m{\ke\ot\ki_X\to\ke} \ induces
\[\BS{\Psi}_\ke:\ke_{|X}\lra\ke_1\ot L^* \ , \]
which is surjective, and we have an exact sequence
\[\Sigma(\ke): \
0\lra\ke_1\lra\ke^{(1)}\lra\ke_{|X}\hfl{\BS{\Psi}_\ke}{}\ke_1\ot L^*\lra 0 \
 . \]
\end{subsub}
Let \m{\Phi:\E\to E}, \m{\Phi':\E'\to E'} be two surjective morphisms of
vector bundles on $X$. We say that $\Phi$ is isomorphic to \m{\Phi'}, or \
\m{\Phi\simeq\Phi'}, if there exists isomorphisms \m{\A:\E\to\E'}, \m{A:E\to
E'} such that the following diagram is commutative
\xmat{\E\ar[rr]^-\Phi\ar[d]^\A & & E\ar[d]^A\\ \E'\ar[rr]^-{\Phi'} & & E' }

Let
\begin{equation}\label{equ14}\Sigma: \ 0 \lra E\ot
L\hfl{\alpha}{}\F\hfl{\beta}{}\E
\hfl{\Phi}{} E\lra 0 \ , \end{equation}
\[\Sigma': \ 0 \lra E'\ot L\hfl{\alpha'}{}\F'\hfl{\beta'}{}\E'
\hfl{\Phi'}{} E'\lra 0 \ , \]
be exact sequences of vector bundles on $X$. We say that $\Sigma$ is
isomorphic to \m{\Sigma'}, or \ \m{\Sigma\simeq\Sigma'}, if there are
isomorphisms \m{\A:\E\to\E'}, \m{A:E\to E'}, \m{\B:\F\to\F'} such that the
following diagram is commutative
\xmat{E\ot L\ar[r]^-\alpha\ar[d]^{A\ot I_L} & \F\ar[r]^-\beta\ar[d]^\B &
\E\ar[r]^-\Phi\ar[d]^\A & E\ar[d]^A\\
E'\ot L\ar[r]^-{\alpha'} & \F'\ar[r]^-{\beta'} & \E'\ar[r]^-{\Phi'} & E'}

So we have two problems:
\begin{enumerate}
\item[(i)] Given a surjective morphism of vector bundles on $X$, \m{\Phi:\E\to
E}, when does there exist a quasi locally free sheaf $\ke$ on \m{\BS{\kx}}
such that \ \m{\BS{\Psi}_\ke\simeq\Phi} \ ?
\item[(ii)] Given an exact sequence $\Sigma$ of vector bundles on $X$ as
above, when does there exist a quasi locally free sheaf $\ke$ on \m{\BS{\kx}}
such that \ \m{\Sigma(\ke)\simeq\Sigma} \ ?
\end{enumerate}

\sepprop

\begin{subsub} The obstructions -- \rm For (i), we define a linear map \
\m{\nabla_{\E,E}:\Hom(\E,E)\to\Ext^2_X(\E,E\ot L)} \ such that the quasi
locally free sheaf $\ke$ exists if and only if \ \m{\nabla_{\E,E}(\Phi)=0}.
\end{subsub}

For (ii) we associate to $\Sigma$ an element \
\m{\nabla''(\Sigma)\in\Ext^2_X(E,E\ot L)} \ such
that there exists a quasi locally free sheaf $\ke$ on \m{\BS{\kx}} such that \
\m{\Sigma\simeq\Sigma(\ke)} \ if and only if \ \m{\nabla''(\Sigma)=0}.

Let \ \m{\mu_\Phi:\Ext^2_X(E,E\ot L)\to\Ext^2_X(\E,E\ot L)} \ be the map
induced by $\Phi$. We have of course, for the exact sequence $(\ref{equ14})$,
\m{\nabla_{\E,E}(\Phi)=\mu_\Phi(\nabla''(\Sigma))}.

\sepprop

\begin{subsub} Vector bundles -- \rm A quasi locally free sheaf $\ke$ on
$\BS{\kx}$ is locally free if and only \m{\BS{\Psi}_\ke} is an isomorphism. In
this case, \ie \m{\E=E} and \m{\phi=I_E}, there is already an invariant in
\m{\Ext^2_X(E,E\ot L)}, \m{\nabla(E)}, constructed in \cite{dr10}, that
vanishes if and only $E$ can be extended to a vector bundle on \m{\BS{\kx}},
cf. \ref{ext_hig}, and we have \ \m{\nabla(E)=\nabla_{\E,\E}(I_\E)}. Let
\[\eta_\Phi:\Ext^2_X(\E,\E\ot L)\lra\Ext^2_X(\E,E\ot L) \ , \qquad
\mu_\Phi:\Ext^2_X(E,E\ot L)\lra\Ext^2_X(\E,E\ot L)\]
be the maps induced by $\Phi$. We will prove (cf. \ref{theo3})
\end{subsub}

\sepprop

{\bf Theorem 1:} We have \
\m{\nabla_{\E,E}(\Phi)=\eta_\Phi(\nabla(\E))=\mu_\Phi(\nabla(E))} .

\sepprop

In particular, if \ \m{\Phi:\E\to E} \ is a surjective morphism of vector
bundles on $X$, and if $\E$ (or $E$) can be extended to a vector bundle on
$\BS{\kx}$, then there exists a quasi locally free sheaf $\ke$ on $\BS{\kx}$
such that \ \m{\BS{\Psi}_\ke\simeq\Phi}. This can be proved directly in a
simple way (cf. \ref{coro3}).

\sepprop

\begin{subsub} Construction of the obstructions -- The case of
\m{\nabla_{\E,E}(\Phi)} -- \rm There are two ways. The first and easy one is
to use {\em the canonical exact sequence of extensions}
\[0\lra\Ext^1_X(\E,E\ot L))\lra\Ext^1_{\BS{\kx}}
(\E,E\ot L)\hfl{\gamma}{}\Hom(\E,E)\hfl{\widetilde{\nabla}}{}\Ext^2_X(\E,E\ot
L)\]
which can be built with \v Cech cohomology, and is also a consequence of the
{\em Ext spectral sequence}, using the fact that we have \
\m{\EExt^1_{\BS{\kx}}(\E,E\ot L)\simeq\HHom(\E,E)} (by theorem
\ref{theo5}).

Let \ \m{\sigma\in\Ext^1_{\BS{\kx}}(\E,E\ot L)}, and \ \m{0\to E\ot
L\to\ke\to\E\to 0} \ the corresponding extension. Then the canonical morphism
\ \m{\ke\ot L\to\ke} \ induces \ \m{\Phi=\gamma(\sigma):\E\to E}, which is
surjective if and only if \ \m{\E=\ke_{|X}}. In this case we have \
\m{\nabla_{\E,E}(\Phi)=\widetilde{\nabla}(\Phi)}.

The second way uses the cohomology theory of sheaves of groups (cf. \cite{fr},
\ref{coh_gr}). Let \ \m{m=\rk(E)}, \m{n=\rk(\E)-m}. To build
\m{\nabla_{\E,E}(\Phi)} we use a sheaf \m{\kg_{m,n}} of algebraic groups on
$X$ such that \m{H^1(X,\kg_{m,n})} is naturally identified with the set of
isomorphism classes of quasi locally free sheaves on \m{\BS{\kx}}, locally
isomorphic to \ \m{(\ko_{\BS{\kx}}\ot\C^m)\oplus(\ko_X\ot\C^n)}. We explain
this method in detail. The same is used in the case of \m{\nabla''(\Sigma)},
and in this case we will simply outline the main points.

For every open subset \m{U\subset
X}, \m{\kg_{m,n}(U)} is the group of automorphisms of \Nligne
\m{(\ko_{U^{(2)}}\ot\C^m)\oplus(\ko_U\ot\C^n)}.

Let $\kh_{m,n}$ be the sheaf
of algebraic groups on $X$, such that for every open subset \m{U\subset X},
\m{\kh_{m,n}(U)} is the group of automorphisms
$\phi$ of the \m{\ko_X}-module \ \m{\ko_X(U)\ot(\C^m\oplus\C^n)} \ such that \
\m{\phi(\ko_X\ot\C^n)\subset\ko_X\ot\C^n}. The cohomology set
\m{H^1(X,\kh_{m,n})} is the set of isomorphism classes of surjective morphisms
\ \m{\E\to E} \ of vector bundles on $X$, where $\E$, $E$ are of rank \m{m+n},
$n$ respectively.

There is a natural surjective morphism of sheaves of algebraic groups,
\m{\rho:\kg_{m,n}\to\kh_{m,n}}. For every \m{g\in H^1(X,\kg_{m,n})},
corresponding to the sheaf $\ke$, \m{h^1(\rho)(g)\in H^1(X,\kh_{m,n})} \
corresponds to \m{\BS{\Psi}_\ke}. The kernel of $\rho$ is the sheaf of additive
groups \ \m{L\ot\Hom(\C^m\oplus\C^n,\C^m)}. Let \ \m{h\in H^1(X,\kh_{m,n})},
\ie a surjective morphism \ \m{\Phi:\E\to E}, where $\E$ (resp. $E$) is a
vector bundle of rank \m{m+n} (resp. $m$) on $X$. From proposition
\ref{prop1b}, \m{\ker(\rho)^h} is isomorphic to the sheaf of commutative
groups \m{\E^*\ot E\ot L}. By the theory of cohomology of sheaves of groups
(cf. \ref{coh_gr}) one associates to $h$ an element \m{\nabla_0} of
\m{H^2(X,\ker(\rho)^h)} which vanishes if and only \ \m{h\in\imm(h^1(\rho))},
\ie if and only if there exists a quasi locally free sheaf $\ke$ such that \
\m{\Phi\simeq\BS{\Psi}_\ke}. We will prove that $\nabla_0$ is the same as
\m{\nabla_{\E,E}(\Phi)} defined by the first way.

\end{subsub}

\sepprop

\begin{subsub} Construction of the obstructions -- The case of
\m{\nabla''(\Sigma)} -- \rm We did not find a simple proof, and use
\m{\kk_{m,n}}, the sheaf of algebraic groups on
$X$, such that for every open subset \m{U\subset X}, \m{\kk_{m,n}(U)} is the
group of triples
\[\begin{pmatrix}
A_0\ot I_L \ , \quad \begin{pmatrix}A_0\ot I_L & B\\ 0 & E\end{pmatrix} \ ,
\quad\begin{pmatrix}A_0 & 0 \\ C & E\end{pmatrix}\end{pmatrix} \ , \]
where \m{A_0} (resp. $E$) is an automorphism of \ \m{\ko_U\ot\C^m} (resp. \
\m{\ko_U\ot\C^n}), \Nligne \m{B\in\Hom(\ko_U\ot\C^n, L\ot\C^m)},
\m{C\in\Hom(\ko_U\ot \C^m,\ko_U\ot\C^n)}. Then \m{H^1(X,\kk_{m,n})} is the set
of isomorphism classes of exact sequences $(\ref{equ14})$. We have a natural
morphism \ \m{\Theta:\kg_{m,n}\to\kk_{m,n}}, and \
\m{\ker(\Theta)=L\ot\End(\C^m)}. Let \ \m{h\in  H^1(X,\kh_{m,n})},
corresponding to $\Sigma$. We have \ \m{\ker(\Theta)^h\simeq E^*\ot E\ot L},
whence the definition of \m{\nabla''(\Sigma)}.

\end{subsub}

\end{sub}

\sepsub

\Ssect{Families of quasi locally free sheaves}{Intro2}

To build such faimiles we use {\em relative {\rm Ext}-sheaves} studied in
\cite{la}.
Let $S$, $T$ be smooth irreducible varieties and $\E$ (resp. $\F$) a vector
bundle on \m{X\times S} (resp. \m{X\times T}).
Let \m{p_S:X\times S\times T\to X\times S}, \m{p_T:X\times S\times T\to X\times
T}, \m{p_{S\times T}:\BS{\kx}\times S\times T\to S\times T} be the
projections. We view \m{p_S^*(\E)} and \m{p_T^*(\F)} as sheaves on \
\m{\BS{\kx}\times S\times T}. We prove that if suitable properties are
satisfied by $\E$ and $\F$, then we have

\sepprop

{\bf Theorem 2:} {\em \m{{\bf Ext}^1_{p_{S\times T}}(p_S^*(\E), p_T^*(\F))}
commutes with base change.}

\sepprop

This can be used to build universal families of quasi locally free sheaves.

In \ref{theo6} we demonstrate what justifies the name ``rigid quasi locally
free sheaf`` :

\sepprop

{\bf Theorem 3: } {\em Let $Y$ be a scheme, $\ke$ a
coherent sheaf on \m{\BS{\kx}\times Y}, flat on $Y$. Let \m{y_0} be a closed
point of $Y$ such that \m{\ke_{y_0}} is quasi locally free of rigid type. Then
there exists an open subset \m{U\subset Y} such that \m{y_0\in U}, and
\begin{enumerate}
\item[--] For every closed point \m{y\in U}, \m{\ke_y} is quasi locally free of
rigid type, of the same type as \m{\ke_{y_0}}.
\item[--] $\ke_{|X\times U_{red}}$ and \m{\ke_{1|X\times U_{red}}} are vector
bundles on \m{X\times U_{red}}, and for every closed point \m{y\in U}
\[(\ke_{|U_{red}\times X})_y \ = \ \ke_{y|X} \ , \quad (\ke_{1|U_{red}\times
X})_y \ =(\ke_y)_1 \ . \]
\end{enumerate}}

\sepprop

Let $R$ be a positive integer. If $R$ is even, the simplest ``generic`` sheaves
on \m{\BS{\kx}} are the vector bundles of rank \m{R/2}. If $R$ is odd, the
simplest generic sheaves on \m{\BS{\kx}} are the rigid quasi locally free
sheaves of type \m{(1,\frac{R-1}{2})}.

Let $C$ be a smooth curve. We can see \ \m{\BS{\kx}\times C} \ in an obvious
way as a primitive double scheme. For families of quasi locally free sheaves
of the same type we have the following result:

\sepprop

{\bf Theorem 4: } {\em Let $\ke$ be a coherent sheaf on \m{\BS{\kx}\times C},
flat on $C$. Suppose that for every closed point \m{c\in C}, \m{\ke_c} is
quasi locally free, and that the type of \m{\ke_c} is independent of $c$. Then

{\bf 1 -- } $\ke$ is quasi locally free.

{\bf 2 -- } The sheaves \m{\ke_1} and \m{\ke_{|X\times C}} are flat on $C$, and
for every \m{c\in C}, the morphism \ \m{(\ke_1)_c\to\ke_c}
\ induced by \ \m{\ke_1\subset\ke} is injective, and
\[(\ke_1)_c \ = \ (\ke_c)_1 \ , \quad \ (\ke_{|X})_c \ =
\ke_{c|X} \ . \]

{\bf 3 -- } The sheaf \m{\ke^{(1)}} is flat on $C$, and for every \m{c\in
C}, the morphism \ \m{(\ke^{(1)})_c\to\ke_c} \ induced by \
\m{\ke^{(1)}\subset\ke} is injective, and \
\m{(\ke^{(1)})_c \ = (\ke_c)^{(1)}}.

{\bf 4 -- } For every \m{c\in C}, the image of the
Koda\"ira-Spencer morphism \
\m{TC_c\to\Ext^1_{\BS{\kx}}(\ke_c,\ke_c)} \ is contained in the canonical
subspace \ \m{H^1(\BS{\kx},\EEnd(\ke_c))} \ of \
\m{\Ext^1_{\BS{\kx}}(\ke_c,\ke_c)}.
}

\sepprop

Statements {\bf 1}, {\bf 2} and {\bf 3} are proved in \cite{dr12}, and {\bf 4}
in \ref{coco5}. In general \m{H^1(\BS{\kx},\EEnd(\ke_c))} is a proper
subspace, hence in this case one cannot construct complete families or moduli
spaces of rigid sheaves by only giving a parameterization that comes from
their construction. There is still a non reduced structure on such complete
families or moduli spaces. When $X$ is a curve, in some cases this non reduced
structure is related to deformations of $\BS{\kx}$ (cf. \cite{dr9}).

\end{sub}

\sepsub

\Ssect{Examples}{Intro3}

We will consider the first non trivial case, quasi locally free sheaves
locally isomorphic to \ \m{\ko_{\BS{\kx}}\oplus\ko_X}, and take for $X$,
\m{\P_2} or \m{\P_1\times\P_1}, and for $L$ the canonical line bundle
\m{\omega_X}.

\sepprop

\begin{subsub} The case \ \m{X=\P_2} -- \rm The are only two possibilities for
$\BS{\kx}$; either it is the trivial double scheme \m{\BS{X}_2^0}, or a non
projective scheme \m{\BS{X}_2}.

Every line bundle \m{\ko_{\P_2}(m)} extends uniquely to a line bundle
\m{\ko_{\BS{X}_2^0}(m)} on \m{\BS{X}^0_2}. On \m{\BS{X}_2}, the only line
bundle is the trivial one.

We prove in \ref{prop10} that the only quasi locally free sheaves locally
isomorphic to \ \m{\ko_{\BS{\kx}}\oplus\ko_X} are
\begin{enumerate}
\item[--] On \m{\BS{X}_2^0} : \m{\ko_{\BS{X}^0_2}(m)\oplus\ko_{\P_2}(n)} ,
\m{m,n\in\Z}.
\item[--] On \m{\BS{X}_2} : \m{\ko_{\BS{X}_2}\oplus\ko_{\P_2}(n)} ,
\m{n\in\Z}.
\end{enumerate}

\end{subsub}

\sepprop

\begin{subsub} The case \ \m{X=\P_1\times\P_1} -- \rm We have \
\m{H^1(X,T_X\ot L)=H^1(X,\Omega_X)\simeq\C^2}, generated by the classes of the
line bundles \m{\ko(1,0)} and \m{\ko(0,1)}. The non trivial primitive double
schemes with reduced subscheme $X$ and associated line bundle \m{\omega_X} are
parameterized by \m{\P(H^1(X,\Omega_X))\simeq\P_1}. Among this set of schemes,
there are projective and non quasi projective ones. Their Picard groups can be
isomorphic to \m{\Z^2}, \m{\Z}, or contain only \m{\ko_X}. If a line bundle on
$X$ can be extended fo $\BS{\kx}$, the extension is unique (cf. \ref{p1p1}).

Let $\ke$ be a quasi locally free sheaf on $\BS{\kx}$, locally isomorphic to \
\m{\ko_{\BS{\kx}}\oplus\ko_X}, and \Nligne
\m{\BS{\Psi}_\ke:\ke_{|X}\to\ke_1\ot\omega_X^*} \ the canonical morphism. We
have \ \m{\ke_1\ot\omega_X^*\simeq\ko(a_1,b_1)} \ and \Nligne
\m{\ker(\BS{\Psi}_\ke)=\ko(a_0,b_0)}, for some \m{a_0,b_0,a_1,b_1\in\Z}, \ie
we have an exact sequence
\[0\lra\ko(a_0,b_0)\lra\ke_{|X}\hfl{\BS{\Psi}_\ke}{}\ko(a_1,b_1)\lra 0 \ . \]
If this extension trivial, then \m{\BS{\Psi}_\ke} is the projection \
\m{p:\ko(a_0,b_0)\oplus\ko(a_1,b_1)\to\ko(a_1,b_1)}. We prove that there
exists a quasi locally free sheaf $\ke$ on \m{\BS{\kx}_2} such that \
\m{\BS{\Psi}_\ke\simeq p} \ if and only if \m{\ko(a_1,b_1)} can be extended to
a line bundle on \m{\BS{\kx}_2}.

Now consider non trivial extensions
\begin{equation}0\lra\ko(a_0,b_0)\lra\E\hfl{\Phi}{}
E=\ko(a_1,b_1)\lra 0 \ . \end{equation}
They exist only if
\begin{equation}
a_0-a_1 \ \geq 0 \ , \ b_0-b_1\leq -2 \qquad \text{or} \qquad
a_0-a_1 \ \leq \ -2 \ , \ b_0-b_1\geq 0 \ .\end{equation}
In this case there always exist quasi locally free sheaves $\ke$ on
\m{\BS{\kx}_2} such that \ \m{\BS{\Psi}_\ke\simeq\Phi}. Let \
\m{\Sigma(a_0,b_0,a_1,b_1)} \ be the set of isomorphism classes of quasi
locally free sheaves $\ke$ on \m{\BS{\kx}_2} such that \
\m{\ke_1\simeq\ko(a_0-2,b_0-2)} \ and \ \m{\ke^{(1)}/\ke_1\simeq\ko(a_1,b_1)}.
From \ref{Intro2}, deformations of sheaves of this type are also of this type.
Suppose that \m{a_0\not=a_1} and \m{b_0\not=b_1}. We build a smooth
irreducible variety \m{\pmb{\M}'}, a coherent sheaf $\kh$ on
\m{\BS{\kx}\times\pmb{\M}'}, flat on \m{\pmb{\M}'}, such that
\begin{enumerate}
\item[--] For every \m{m\in\pmb{\M}'}, \m{\kh_m} is a sheaf of
\m{\Sigma(a_0,b_0,a_1,b_1)}, and the Koda\"ira-Spencer map \
\m{T_{\pmb{\M}',m}\to\Ext^1_{\BS{\kx}_2}(\kh_m,\kh_m)} \ is an isomorphism \
\m{T_{\pmb{\M}',m}\simeq H^1(\BS{\kx}_2,\EEnd(\kh_m))}.
\item[--] For every integral scheme $S$ and every flat family $\ke$ of
sheaves in \m{\Sigma(a_0,b_0,a_1,b_1)} parameterized by $S$, there exists a
morphism \ \m{f_\ke:S\to\pmb{\M}'} \ such that for every closed point \m{s\in
S}, there exists a neighborhood $U$ of $s$ and an isomorphism \
\m{\ke_{|\BS{\kx}\times U}\simeq f_\ke^*(\kh)_{|\BS{\kx}\times U}}.
\end{enumerate}
We call this a fine moduli space for \m{\Sigma(a_0,b_0,a_1,b_1)}, for families
of sheaves parameterized by integral schemes. But we have \
\m{H^1(\BS{\kx}_2,\EEnd(\kh_m))\not=\Ext^1_\BS{\kx}(\kh_m,\kh_m)}. So a true
fine moduli space, if it exists, is a non reduced structure on \m{\pmb{\M}'}.

\end{subsub}

\end{sub}

\sepsub

\Ssect{Outline of the paper}{Intro4}

In Chapter 2, we give several preliminary definitions and recall some results
that will be used in the next chapters.

In Chapter 3 we recall the construction of primitive double schemes and some
of their properties.

Chapter 4 is devoted to the study of quasi locally free sheaves, and to the
proof of theorems 3 and 4.

In Chapter 5 we define the map \m{\nabla_{\E,E}} which solves problem (i) of
\ref{Intro1}, using the first way (the canonical exact sequence of
extensions), and prove theorem 1.

In Chapter 6 we define again the map \m{\nabla_{\E,E}}, and also
\m{\Sigma(\ke)}, using the cohomology theory of sheaves of groups, which
solves problem (ii).

Chapter 7 and Chapter 8 are devoted to the examples described in \ref{Intro3}.

\end{sub}

\sepsec

\section{Preliminaries}\label{prelim}

\Ssect{\v Cech cohomology and trivializations}{cech}

Let $X$ be a scheme over $\C$ and $E$ a vector bundle on $X$, \m{r=\rk(E)}.
Let \m{(U_i)_{i\in I}} be an open cover of $X$ such that we have
trivializations:
\xmat{\alpha_i:E_{|U_i}\ar[r]^-\simeq & \ko_{U_i}\ot\C^r .}
Let
\xmat{\alpha_{ij}=\alpha_i\circ\alpha_j^{-1}:\ko_{U_{ij}}\ot\C^r\ar[r]^-\simeq &
\ko_{U_{ij}}\ot\C^r ,}
so that we have the relation \ \m{\alpha_{ij}\alpha_{jk}=\alpha_{ik}}.

\sepprop

\begin{subsub}\label{cech1}\rm Let $n$ be a positive integer, and for every
sequence \m{(i_0,\ldots,i_n)} of distinct elements of $I$, \m{\sigma_{i_0\cdots
i_n}\in H^0(U_{i_0\cdots i_n},\ko_X\ot\C^r)}. Let
\[\theta_{i_0\cdots i_n} \ = \ \alpha_{i_0}^{-1}\sigma_{i_0\cdots i_n} \ \in \
H^0(U_{i_0\cdots i_n},E) \ . \]
The family \m{(\theta_{i_0\cdots i_n})} represents an element of \m{H^n(X,E)}
if the cocycle relations are satisfied: for every sequence
\m{(i_0,\ldots,i_{n+1})} of distinct elements of $I$,
\[\sigg_{k=0}^n(-1)^k\theta_{i_0\cdots\widehat{i_k}\cdots i_{n+1}} \ = \ 0 \ ,
\]
which is equivalent to
\[\alpha_{i_0i_1}\sigma_{i_1\cdots i_{n+1}}+\sigg_{k=1}^{n+1}
(-1)^k\sigma_{i_0\cdots\widehat{i_k}\cdots i_{n+1}} \ = \ 0 \ . \]
For \m{n=1}, this gives that elements of \m{H^1(X,E)} are represented by
families \m{(\sigma_{ij})},\Nligne \m{\sigma_{ij}\in H^0(U_{ij},\ko_X\ot\C^r)},
such that
\[\alpha_{ij}\sigma_{jk}+\sigma_{ij}-\sigma_{ik} \ = \ 0 \ . \]
In \v Cech cohomology, it is generally assumed that
\m{\theta_{ji}=-\theta_{ij}}. This implies that \
\m{\sigma_{ji}=-\alpha_{ji}\sigma_{ij}}.
\end{subsub}

\sepprop

\begin{subsub}\label{cech2}\rm Let $F$ be a vector bundle on $X$. Similarly,
an
element of \m{H^n(X,F\ot E)} is represented by a family  \m{(\mu_{i_0\cdots
i_n})}, with \ \m{\mu_{i_0\cdots i_n}\in H^0(U_{i_0\cdots i_n},F\ot\C^r)},
satisfying the relations
\[(I_F\ot\alpha_{i_0i_1})(\mu_{i_1\cdots i_{n+1}})+\sigg_{k=1}^{n+1}
(-1)^k\mu_{i_0\cdots\widehat{i_k}\cdots i_{n+1}} \ = \ 0 \ . \]
The corresponding element of \m{H^0(U_{i_0\cdots i_n},E\ot F)} is \
\m{\theta_{i_0\cdots i_n}=(I_F\ot\alpha_{i_0})^{-1}(\mu_{i_0\cdots i_n})}.

For \m{n=1}, this gives that elements of \m{H^1(X,F\ot E)} are represented by
families \m{(\sigma_{ij})},\Nligne \m{\sigma_{ij}\in H^0(U_{ij},F\ot\C^r)},
such that
\[(I_F\ot\alpha_{ij})\sigma_{jk}+\sigma_{ij}-\sigma_{ik} \ = \ 0 \ . \]
If $E$ is a line bundle, then \ \m{\alpha_{ij}\in\C^*}, and the cocycle
relation is \ \m{\alpha_{ij}\sigma_{jk}+\sigma_{ij}-\sigma_{ik}=0}.

The family \m{(\sigma_{ij})} represents 0 if and only if there exists, for
every \m{i\in I}, \m{\rho_i\in H^0(U_i,F\ot\C^r)}, such that \
\m{\sigma_{ij}=\rho_i-(I_F\ot\alpha_{ij})(\rho_j)}.
\end{subsub}

\end{sub}

\sepsub

\Ssect{Extensions of coherent sheaves}{ext_sh}

Let $U$ be an affine scheme and $E$, $F$ coherent sheaves on $U$. Suppose that
we have a locally free resolution of $E$:
\xmat{\cdots\E_2\ar[r]^-{\epsilon_2} & \E_1\ar[r]^-{\epsilon_1} &
\E_0\ar[r]^-{\epsilon_0} & E \ . }
Then we have a canonical isomorphism \ \m{\Ext^1_X(E,F)\simeq H/K}, where
\ \m{H,K\subset\Hom(\E_1,F)}, $H$ is the space of morphisms that vanish on
\m{\imm(\epsilon_2)}, and $K$ is the space of morphisms of the form
\m{\phi\circ\epsilon_1}, for some \ \m{\phi:\E_0\to F}.

Let \ \m{\sigma\in\Ext^1_X(E,F)}, defined by \ \m{\alpha\in H}, and
\[0\lra F\hfl{a}{}\ke\hfl{b}{} E\lra 0\]
the corresponding extension. The sheaf $\ke$ is defined as follows: let
\[\mu=(\alpha,\epsilon_1):\E_1\lra F\oplus\E_0 \ . \]
Then \ \m{\ke=\coker(\mu)}, $a$ and $b$ being defined in an obvious way.

Let $Y$ be a noetherian scheme, and
\[0\lra A\lra\ke\lra B\lra 0\]
and exact sequence of coherent sheaves on $Y$, corresponding to \
\m{\sigma\in\Ext^1_Y(B,A)}. Let \m{C} be a coherent sheaf on $Y$ and \
\m{f:C\to B} \ a morphism, and \
\m{\phi:\Ext^1_Y(B,A)\to\Ext^1_Y(C,A)} \ the induced map. Let
\[0\lra A\lra\kf\lra C\lra 0\]
an exact sequence of coherent sheaves on $Y$, corresponding to \
\m{\tau\in\Ext^1_Y(C,A)}. Then we have \ \m{\tau=\phi(\sigma)} \ if and
only if there exists a commutative diagram
\xmat{0\ar[r] & A\ar[r]\fleq[d] & \kf\ar[r]\ar[d] & C\ar[r]\ar[d]^f & 0\\
0\ar[r] & A\ar[r] & \ke\ar[r] & B\ar[r] & 0}

\end{sub}

\sepsub

\Ssect{Flat families of sheaves}{FLA}

We give here some results that will be used later.

\sepprop

\begin{subsub}\label{prop6}{\bf Proposition: } Let $X$, $Y$ be schemes and
$\ke$ a coherent sheaf on \m{X\times Y}, flat over $Y$. Then

{\bf 1 -- } Let \m{Y'\subset Y} be a closed subscheme. Then \m{\ke_{|X\times
Y'}} is flat over \m{Y'}.

{\bf 2 -- } Let \m{X'\subset X} be a closed subscheme. Suppose that the support
of $\ke$ is contained in \m{X'\times Y}. Then $\ke$, as a sheaf on \m{X'\times
Y}, is flat over $Y$.

{\bf 3 -- } Let \m{x\in X}, \m{y\in Y} be closed points. Suppose that \m{\ke_y}
is locally free at $x$. Then $\ke$ is locally free in a neighborhood of
\m{(x,y)}.
\end{subsub}
\begin{proof}
{\bf 1} follows from \cite{sta}, lemma 10.39.7.

{\bf 2} follows from the following easy lemma of algebra: let $A$ be a
commutative ring and $J$ an ideal of $A$. Let $M$ be a flat $A$-module such
that \ \m{JM=0}. Then $M$ is a flat \m{A/J}-module.

{\bf 3} follows from lemma 37.16.7 of \cite{sta}.
\end{proof}

\end{sub}

\sepsub

\Ssect{Cohomology of sheaves of non abelian groups}{coh_gr}

(cf. \cite{fr})

Let $X$ be a paracompact topological space and $G$ a sheaf of groups on $X$. We
denote by \m{H^1(X,G)} the first \v Cech cohomology set of $G$. If $G$ is not
a sheaf of {\em commutative} groups there is no natural structure of group on
\m{H^1(X,G)}. Let \m{e\in H^1(X,G)} be defined by the trivial section of $G$
(on the trivial open cover \m{(X)} of $X$), we call $e$ the {\em identity
element} of \m{H^1(X,G)}.

Given an open cover \m{(U_i)_{i\in I}} of $X$, a 1-cocycle of $G$ is a family
\m{(g_{ij})_{i,j\in I}} such that \m{g_{ij}\in G(U_{ij})} and that the
relations \m{g_{ij}g_{jk}=g_{ik}} are satisfied.

\sepsubsub

\begin{subsub}\label{coh_gr1} Actions of sheaves of groups -- \rm
Let $H$ be another sheaf of group on $X$, acting on $G$ (the action is
compatible with the group actions).
Let \m{z\in H^1(X,H)}, represented by a cocycle \m{(z_{ij})}
with respect to an open cover \m{(U_i)_{i\in I}} of $X$. Recall that the
associated sheaf of groups \m{G^z} is defined as follows: we glue the sheaves
\m{G_{|U_i}} using the automorphisms \m{(z_{ij})} of  \m{G_{|U_{ij}}}. The
sections of \m{G^z} are represented by families \m{(\gamma_i)_{i\in I}},
\m{\gamma_i\in G(U_i)}, such that \m{\gamma_i=z_{ij}\gamma_j} on \m{U_{ij}}.
The elements of \m{H^1(X,G^z)} are represented by families \m{(\rho_{ij})},
\m{\rho_{ij}\in G(U_{ij})}, such that \
\m{\rho_{ij}^{-1}\rho_{ik}=z_{ij}\rho_{jk}} \ on \m{U_{ijk}}.
If $G$ is a sheaf of {\em commutative} groups, the elements of the higher
cohomology groups \m{H^k(X,G^z)}, \m{k\geq 2} can be represented in the same
way by appropriate families.
\end{subsub}

\sepsubsub

\begin{subsub}\label{coh_gr2}Exact cohomology sequence -- \rm Let
\m{\Gamma\subset G} be a subsheaf of normal subgroups of $G$, and \m{G/\Gamma}
the quotient sheaf of groups. Hence we have an exact sequence of sheaves of
groups
\xmat{0\ar[r] & \Gamma\ar[r]^-i & G\ar[r]^-p & G/\Gamma\ar[r] & 0 . }

Let \m{A_1,\ldots,A_m} be sets, having a particular element \m{e_i\in A_i} for
\m{1\leq i\leq m}. A sequence of maps
\xmat{A_1\ar[r]^-{f_1} & A_2\ar[r] & \cdots\ar[r] & A_{m-1}
\ar[r]^-{f_{m-1}} & A_m}
is called {\em exact} if for \m{1<i<m} we have \
\m{f_i^{-1}(e_i)=f_{i-1}(A_{i-1})} . Then we have: {\em the canonical sequence
\xmat{0\ar[r] & H^0(X,\Gamma)\ar[r]^-{H^0(i)} & H^0(X,G)\ar[r]^-{H^0(p)}
& H^0(X,G/\Gamma)\ar[r]^-\delta & \null\quad\quad\quad\quad}
\xmat{\quad\quad\quad\quad\quad\quad H^1(X,\Gamma)\ar[r]^-{H^1(i)} &
H^1(X,G)\ar[r]^-{H^1(p)} & H^1(X,G/\Gamma)\quad .}
is exact. The subsequence
\xmat{0\ar[r] & H^0(X,\Gamma)\ar[r]^-{H^0(i)} & H^0(X,G)\ar[r]^-{H^0(p)}
& H^0(X,G/\Gamma)}
is an exact sequence of groups. Moreover, if \m{c,c'\in H^0(X,G/\Gamma)}, we
have \ \m{\delta(c)=\delta(c')} \ if and only if \
\m{c{c'}^{-1}\in\imm(H^0(p))}.}

The group \m{H^0(X,G/\Gamma)} acts on \m{H^1(X,\Gamma)} as follows: let \m{a\in
H^1(X,\Gamma)}, represented by a cocycle \m{(\gamma_{ij})} (with respect to an
open cover \m{(U_i)} of $X$), and \m{c\in H^0(X,G/\Gamma)}. Suppose also that
the cover \m{(U_i)} is such that \m{c_{|U_i}} is the image of an element
\m{c_i} of \m{G(U_i)}. Then \m{(c_i\gamma_{ij}c_j^{-1})} is a cocycle of
elements of $\Gamma$, and the corresponding element of \m{H^1(X,\Gamma)}
depends only on $c$ and $a$; we denote it by \m{c.a}. We have then \
\m{\delta(c)=c^{-1}.e} \ . Moreover, let \ \m{\gamma,\gamma'\in H^1(X,\Gamma)}.
Then we have: {\em\m{H^1(i)(\gamma)=H^1(i)(\gamma')} \ if and only if
\m{\gamma,\gamma'} are in the same \m{H^0(X,G/\Gamma)}-orbit.}
\end{subsub}

\sepsubsub

\begin{subsub}\label{coh_gr3}The fibers of \m{H^1(p)} -- \rm The sheaf of
groups $G$ acts (by conjugation) on itself and $\Gamma$, \m{G/\Gamma}. Let \
\m{\omega\in H^1(X,G/\Gamma)}, \m{g\in H^1(X,G)} be such that \
\m{H^1(p)(g)=\omega}. Suppose that $g$ is represented by a cocycle \m{(g_{ij})}.
Let \m{g'\in H^1(p)^{-1}(\omega)}. Then \m{g'} is represented by a cocycle of
the form \m{(\gamma_{ij}g_{ij})}, with \m{\gamma_{ij}\in\Gamma(U_{ij})}.
Conversely, a cochain \m{(\gamma_{ij}g_{ij})} is a cocycle if and only if, for
every \m{i,j,k} we have \ \m{\gamma_{ij}^{-1}\gamma_{ik} =
g_{ij}\gamma_{jk}g_{ij}^{-1}}, \ie if and only if \m{(\gamma_{ij})}
induces a cocycle in \m{\Gamma^g} (cf. \ref{coh_gr1}). In this way we define a
surjective map
\[\lambda_g:H^1(X,\Gamma^g)\lra H^1(p)^{-1}(\omega)\]
sending the identity element of \m{H^1(X,\Gamma^g)} to $g$.
\end{subsub}
Note that we have an exact sequence \ \m{0\to\Gamma^g\to G^g\to G^g/\Gamma^g\to
0}.

The group \m{H^0(X,G^g/\Gamma^g)} acts on \m{H^1(X,\Gamma^g)} and the fibers of
\m{\lambda_g} are the fibers of this action.

\sepsubsub

\begin{subsub}\label{coh_gr4} The case where $\Gamma$ is commutative -- \rm
The sheaf of groups \m{G/\Gamma} acts then by conjugation on $\Gamma$. Let \
\m{\omega\in H^1(X,G/\Gamma)}. We define \ \m{\Delta(\omega)\in
H^2(X,\Gamma^\omega)} \ as follows: we can represent $\omega$ by a cocycle
\m{(\omega_{ij})} (for a suitable open cover \m{(U_i)} of $X$), such that for
every \m{i,j} there exists \m{g_{ij}\in G(U_{ij})} over \m{\omega_{ij}}, and
that \ \m{g_{ji}=g_{ij}^{-1}}. For every indices \m{i,j,k} let \
\m{\gamma_{ijk}=g_{ij}g_{jk}g_{ki}}. Then \m{(\gamma_{ijk})} is a family
representing \m{\Delta(\omega)} (cf. \ref{coh_gr1}). We have then: {\em
\m{\Delta(\omega)=0} if and only if $\omega$ belongs to the image of
\m{H^1(p)}}.

The commutativity of $\Gamma$ brings also a supplementary property of the map
\m{\lambda_g} of \ref{coh_gr3}. The sheaf of groups \m{\Gamma^g} can be seen
naturally as a subsheaf of commutative groups of of \m{G^g}, and we have a
canonical isomorphism \ \m{(G/\Gamma)^g\simeq G^g/\Gamma^g}. Hence there is a
natural action of \m{H^0(X,G^g/\Gamma^g)} on \m{H^1(X,\Gamma^g)}, and we have:
{\em the fibers of \m{\lambda_g} are the orbits of the action of
\m{H^0(X,G^g/\Gamma^g)}.} This action is as follows: let \m{h\in
H^0(G^g/\Gamma^g)}, \m{\gamma\in H^1(X,\Gamma^g)}, represented by
\m{(\gamma_{ij})}. For a suitable open cover \m{(U_i)}, each \m{h_{|U_i}} can
be lifted to \m{h_i\in G(U_i)}, and we have \
\m{h_ig_{ij}h_jg_{ij}^{-1}\in\Gamma(U_{ij})}. Then \m{h.\gamma} is represented
by \m{(h_i.\gamma_{ij}.g_{ij}h_j^{-1}g_{ij}^{-1})}.
\end{subsub}

\end{sub}

\sepsub

\Ssect{Relative Ext-sheaves}{rel_e}

(cf. \cite{la})

Let $X$ be a proper scheme, $Y$ a scheme, \m{p_Y:X\times Y\to Y} the
projection, and $\kf$, $\kg$ coherent sheaves on \m{X\times Y}, flat on $Y$,
$i$ an integer, \m{i\geq 0}. Let \m{{\bf Ext}^i_{p_Y}(\kf,\kg)} be the sheaf on
$Y$
associated to the presheaf
\[U\mapsto\Ext^i_{X\times U}(\kf_{|X\times U},\kg_{|X\times U}) \ . \]
Then
\begin{enumerate}
\item[(i)] if \m{Y=\spec(A)} is affine, then there is a natural isomorphism
\[{\bf Ext}^i_{p_Y}(\kf,\kg) \ \simeq \ \Ext^i_{X\times Y}(\kf,\kg)^{\sim} \ .
\]
\item[(ii)] if $\kf$ is locally free, then there is a natural isomorphism
\[{\bf Ext}^i_{p_Y}(\kf,\kg) \ \simeq \ R^ip_{Y*}(\kf^*\ot\kg) \ . \]
\end{enumerate}
Let \m{y\in Y} be a closed point. Then there is a canonical morphism
\[\tau^i(y):{\bf Ext}^i_{p_Y}(\kf,\kg)_y\lra\Ext^i_X(\kf_y,\kg_y) \ , \]
which induces, for every open subset \m{U\subset Y} containing $y$, a canonical
map
\[\Phi^i_{U,y}:\Ext^i_{X\times U}(\kf_{|X\times U},\kg_{|X\times U})\lra
\Ext^i_X(\kf_y,\kg_y) \ . \]

Recall theorem 1.4 of \cite{la}:

\sepprop

\begin{subsub}\label{theo8}{\bf Theorem: } Let \m{y\in Y} be a point and assume
that \m{\tau^i(y)} is surjective. Then:
\begin{enumerate}
\item[--] there is a neighborhood $U$ of $y$ such that \m{\tau^i(y')} is an
isomorphism for all \m{y'\in U};
\item[--] the map $\tau^{i-1}(y)$ is surjective if and only if \m{{\bf
Ext}^i_{p_Y}(\kf,\kg)} is locally free in a neighborhood of $y$.
\end{enumerate}
\end{subsub}

\sepprop

If \m{\tau^i(y)} is an isomorphism for every \m{y\in Y}, we say that {\em
\m{{\bf Ext}^i_{p_Y}(\kf,\kg)} commutes with base change}. If this property is
satisfied, several interesting constructions, like {\em universal extensions},
are made in \cite{la}.

From (ii) we have

\sepprop

\begin{subsub}\label{prop7}{\bf Proposition: } Suppose \
\m{r=\dim(\Ext^i_X(\kf_y,\kg_y))} \ does not depend on \m{y\in Y}. Then:
\item[(1)] If \m{{\bf Ext}^i_{p_Y}(\kf,\kg)} commutes with base change, then it
is a locally free sheaf of rank $r$.
\item[(2)] If $\kf$ is locally free then  \m{{\bf Ext}^i_{p_Y}(\kf,\kg)}
commutes with base change.
\end{subsub}

\sepprop

Suppose that the hypotheses of proposition \ref{prop7} are satisfied, and that
\m{{\bf Ext}^i_{p_Y}(\kf,\kg)} commutes with base change. Let \
\m{\alpha:Y'\to Y} \ be a morphism of schemes,
\m{\ov{\alpha}=I_X\times\alpha:X\times
Y'\to X\times Y}, and \ \m{p_{Y'}:X\times Y'\to Y'} the projection. Then \
\m{{\bf Ext}^i_{p_{Y'}}(\ov{\alpha}^*(\kf),\ov{\alpha}^*(\kg))} \ commutes with
base change, and there is a canonical morphism of vector bundles on \m{Y'}
\[\tau^i(\alpha):\alpha^*[{\bf Ext}^i_{p_Y}(\kf,\kg)]\lra{\bf Ext}^i_{p_{Y'}}(
\ov{\alpha}^*(\kf),\ov{\alpha}^*(\kg)) \ . \]
For every \m{y'\in Y'}, if \m{y=\alpha(y')}, then
\[\alpha^*[{\bf Ext}^i_{p_Y}(\kf,\kg)]_{y'} \ = \ {\bf Ext}^i_{p_{Y'}}(
\ov{\alpha}^*(\kf),\ov{\alpha}^*(\kg))_{y'} \ = \
\Ext^i_X(\kf_y,\kg_y) \ , \]
and \m{\tau^i(y')} is the identity.

A {\em family of extensions of $\kf$ by $\kg$ over $Y$} is a family
\m{(e_y)_{y\in Y}}, where \ \m{e_y\in\Ext^1_X(\kf_y,\kg_y)}, such that there is
an open covering \m{(Y_i)_{i\in I}} of $Y$ and for every \m{i\in I} an element
\Nligne \m{\sigma_i\in\Ext^1_{X\times Y_i}(\kf_{|X\times Y_i},\kg_{|X\times
Y_i})} \ such that \ \m{e_y=\Phi_{Y_i,y}^1(\sigma_i)}

\end{sub}

\Ssect{Fine moduli spaces of sheaves}{fin-mod}

Let $\chi$ be a set of isomorphism classes of coherent sheaves on a scheme $Z$.
Suppose that $\chi$ is {\em open}, \ie for every scheme $V$ and coherent sheaf
$\E$ on \m{Z\times V}, flat on $V$, if \m{v\in V} is a closed point such that
\m{\E_v\in\chi}, then there exists an open neighborhood $U$ of $v$ such that
\m{\E_u\in\chi} for every closed point \m{u\in U}.
A {\em fine moduli space for $\chi$ }is the data of a scheme $M$ and of
\begin{enumerate}
	\item[--] a bijection
	\[\xymatrix@R=5pt{M^0\ar[r] & \chi\\ m\fmaps[r] & E_m
	}\]
	(where \m{M^0} denotes the set of closed points of $M$),
	\item[--] an open cover \m{(M_i)_{i\in I}} of $M$, and for every \m{i\in I}, a
	coherent sheaf \m{\ke_i} on \m{Z\times M_i} such that for every \m{m\in M_i},
	\m{\ke_{i,m}\simeq E_m},
\end{enumerate}
such that:
for any scheme $S$, any coherent sheaf $\kf$ on \m{Z\times S}, flat on $S$,
such that for any closed point \m{s\in S}, \m{\kf_s\in\chi}, there is a
morphism \ \m{f_\kf:S\to M} \ such that for every \m{s\in S}, if \m{m=f_\kf(s)}
then \m{\kf_s\simeq E_m}, and if \m{m\in M_i}, then there exists an open
neighborhood \m{U} of $s$ such that \ \m{f_\kf(U)\subset M_i} and
\ \m{(I_Z\times f_{\kf|U})^*(\ke_i)\simeq\kf_{|Z\times U}}.

We extend this definition to {\em fine moduli space for families of sheaves in $\chi$ parameterized by smooth schemes} (resp. {\em by integral schemes}), where we consider only smooth (resp. integral) $S$.

\end{sub}

\sepsec

\section{Primitive multiple schemes}\label{PMS}

\Ssect{Definition and construction}{PMS_def}

Let $X$ be a smooth connected variety, and \ \m{d=\dim(X)}. A {\em multiple
scheme with support $X$} is a Cohen-Macaulay scheme $\BS{\kx}$ such that
\m{\BS{\kx}_{red}=X}. If $\BS{\kx}$ is quasi-projective we say that it is a
{\em multiple variety with support $X$}.

We say that $\BS{\kx}$ is {\em primitive} if, for every closed point $x$ of
$X$, there exists a smooth variety $S$ of dimension \m{d+1}, containing a
neighborhood of $x$ in $\BS{\kx}$ as a locally closed subscheme. We say that
$\BS{\kx}$ is a {\em primitive double scheme} if it is primitive,
\m{\ki_X\not=0} and \m{\ki_X^2=0}. In this case \m{L=\ki_X} is a line bundle
on $X$. We call it the line bundle on $X$ {\em associated} to $\BS{\kx}$.

Let \m{P\in X}. Then there exist elements \m{y_1,\ldots,y_d}, $t$ of
\m{m_{S,P}} whose images in \m{m_{S,P}/m_{S,P}^2} form a basis, and such that
for \m{1\leq i<n} we have \ \m{\ki_{X,P}=(t)}. In this case the images of
\m{y_1,\ldots,y_d} in \m{m_{X,P}/m_{X,P}^2} form a basis of this vector space.

Even if $X$ is projective, we do not assume that $\BS{\kx}$ is projective.

The simplest case is when $\BS{\kx}$ is contained in a smooth variety $S$ of
dimension \m{d+1}. Let \m{P\in X} and \m{f\in\ko_{S,P}}  a local equation of
$X$. Then if \m{\ki_{\BS{\kx},S}} is the ideal sheaf of \m{\BS{\kx}} in S, we
have \ \m{\ki_{\BS{\kx},S,P}=(f^2)}, and \ \m{L=\ko_X(-X)} .

For any \m{L\in\Pic(X)}, the {\em trivial double scheme} with induced smooth
variety $X$ and associated line bundle $L$ on $X$ is the second infinitesimal
neighborhood of $X$, embedded by the zero section in the dual bundle $L^*$,
seen as a smooth variety. In this case there is a canonical morphism
\m{\BS{\kx}\to X} which is the identity on $X$.

If \m{U\subset X} is open, we note \m{U^{(2)}} the corresponding open subscheme
of \m{\BS{\kx}}. More generally, if $Y$ is a scheme, and $V$ an open subset of
\ \m{X\times Y}, let \m{V^{(2)}} denotes the corresponding open subscheme
of \ \m{\BS{\kx}\times Y}

\sepsubsub

\begin{subsub}\label{PMS-1} Construction of primitive multiple schemes -- \rm
Let \ \m{{\bf Z}_2=\spec(\C[t]/(t^2))}. If $Y$ is an affine scheme,
\m{Y=\spec(A)}, then \ \m{Y\times{\bf Z}_2=\spec(A[t]/(t^2))}, and every \
\m{\alpha\in\Gamma(Y\times{\bf Z}_2)} \ can be written in a unique way in the
form \ \m{\alpha_0+t\alpha_1}, with \ \m{\alpha_0,\alpha_1\in\Gamma(Y)}.

Let $\BS{\kx}$ be a primitive double scheme, \m{X=\BS{\kx}_{red}} and \
\m{L=\ki_{X,\BS{\kx}}}. Then for every
closed point \m{P\in X}, there exists an open neighborhood $U$ of $P$ in $X$
such that
\begin{enumerate}
\item[--] There exists a section \ \m{\ko_X(U)\to\ko_\BS{\kx}(U^{(2)})} \ of
the restriction map \ \m{\ko_\BS{\kx}(U^{(2)})\to\ko_X(U)},
\item[--] $L_{|U}$ is trivial.
\end{enumerate}
Then there exists a commutative diagram
 \xmat{ & U\flinc[ld]\flinc[rd] \\
U^{(2)}\ar[rr]^-\simeq & & U\times {\mathbf Z}_2}
\ie $\BS{\kx}$ is locally trivial (\cite{dr1}, th\'eor\`eme 5.2.1, corollaire
5.2.2).

It follows that we can construct a primitive multiple scheme of multiplicity
$n$ by taking an open cover \m{(U_i)_{i\in I}} of $X$ and gluing the varieties
\ \m{U_i\times{\bf Z}_2} (with automorphisms of the \ \m{U_{ij}\times{\bf Z}_2}
\ leaving \m{U_{ij}} invariant).

Let \m{(U_i)_{i\in I}} be an affine open cover of $X$ such that we have
trivializations
\xmat{\delta_i:U_i^{(2)}\ar[r]^-\simeq & U_i\times{\bf Z}_2 , }
and \ \m{\delta_i^*:\ko_{U_i\times{\bf Z}_2}\to\ko_{U_i^{(2)}}} \ the
corresponding isomorphism. Let
\xmat{\delta_{ij}=\delta_j\delta_i^{-1}:U_{ij}\times{\bf Z}_2\ar[r]^-\simeq &
U_{ij}\times{\bf Z}_2 \ . }
Then \ \m{\delta_{ij}^*=\delta_i^{*-1}\delta_j^*} \ is an automorphism of \
\m{\ko_{U_i\times Z_2}=\ko_X(U_{ij})[t]/(t^2)}, such that for every \
\m{\phi\in\ko_X(U_{ij})[t]/(t^2)}, seen as a polynomial in $t$ with
coefficients in \m{\ko_X(U_{ij})}, the term of degree zero of
\m{\delta_{ij}^*(\phi)} is the same as the term of degree zero of $\phi$.

If $X$ is affine, then every primitive multiple scheme $\BS{\kx}$ such that \
\m{X=\BS{\kx}_{red}} \ is trivial, \ie there exists an isomorphism \
\m{\BS{\kx}_{red}\simeq X\times{\bf Z}_2} \ inducing the identity on $X$.
\end{subsub}

\sepsubsub

\begin{subsub}\label{I_X} The ideal sheaf of $X$ -- \rm
There exists \ \m{\alpha_{ij}\in\ko_X(U_{ij})} \ such that
\ \m{\delta_{ij}^*(t)=\alpha_{ij}t}.
For every \m{i\in I}, \m{\delta_i^*(t)} is a generator of
\m{\ki_{X,\BS{\kx}|{U^{(2)}}}}. So we have local trivializations
\[\xymatrix@R=5pt{\lambda_i:\ki_{X,\BS{\kx}|{U_i}}\ar[r] &
\ko_{U_i}\\
\delta_i^*(t)\fmaps[r] & 1}\]
Hence \ \m{\lambda_{ij}=\lambda_i\lambda_j^{-1}:
\ko_{U_{ij}^{(n-1)}}\to\ko_{U_{ij}^{(n-1)}}} \ is the multiplication by
\m{\delta_j^*(\alpha_{ij})}. It follows that \m{(\delta_j^*(\alpha_{ij}))} is
a cocycle representing the line bundle \m{L} on $X$.
\end{subsub}

Let \ \m{\alpha_i:L_{|U_i}\to\ko_{U_i}} \ be
isomorphisms such that \ \m{\alpha_{ij}=\alpha_i\circ\alpha_j^{-1}} \ on
\m{U_{ij}}. Then we have the following description of \m{\delta_{ij}^*}:
\begin{enumerate}
\item[--] there is a derivation \m{D_{ij}} of \m{\ko_X(U_{ij})} such that \
\m{\delta_{ij}^*(\beta)=\beta+D_{ij}(\beta)t} \ for every
\m{\beta\in\ko_X(U_{ij})}.
\item[--] $\delta_{ij}^*(t)=\alpha_{ij}t$ .
\end{enumerate}
The relation \ \m{\delta_{ij}^*\delta_{jk}^*=\delta_{ik}^*} \ is equivalent to
\m{D_{ij}+\alpha_{ij}D_{jk}=D_{ik}}, \ie \m{(D_{ij})} is a cocycle in the
sense of \ref{cech2}. We can view \m{D_{ij}} as a section of
\m{T_{X|U_{ij}}}, the family \m{(\alpha_i^{-1}\ot D_{ij})} represents an
element $\lambda$ of \m{H^1(X,T_X\ot L)}, and \m{\C\lambda} is independent of
the choice of the \m{\delta_{ij}^*} and \m{\alpha_i}. We will note \
\m{\C\lambda=\sigma(\BS{\kx})\in\P(H^1(X,T_X\ot L))\cup\nsp}. We have a
canonical exact sequence on $X$
\[ 0\lra L=\ki_{X,\BS{\kx}}\lra\Omega_{\BS{\kx}|X}\lra\Omega_X\lra 0 \ , \]
corresponding to \ \m{\mu\in\Ext^1_X(\Omega_X,L)=H^1(X,T_X\ot L)}. Then we
have \ \m{\sigma(\BS{\kx})=\C\mu}.

\sepprop

\begin{subsub}\label{param_DS} Parameterization of double schemes -- \rm
Let $\BS{\kx}$, \m{\BS{\kx}'} primitive double schemes,
\m{\phi:\BS{\kx}_{red}\to X}, \m{\phi':\BS{\kx}'_{red}\to X} \ isomorphisms. We
say that $\BS{\kx}$, \m{\BS{\kx}'} are are {\em isomorphic} if there is an
isomorphism \ \m{\Phi:\BS{\kx}\to\BS{\kx}'} \ such that \
\m{\phi=\Phi_{red}\circ\phi'}. In this case the associated line bundles on $X$
of $\BS{\kx}$, \m{\BS{\kx}'} are isomorphic.

According to \cite{ba_ei} and \cite{dr10}, two primitive double
schemes \m{X_2}, \m{X'_2}, with underlying smooth variety $X$ and associated
line bundle $L$ are isomorphic (over $X$) if and only if
\m{\sigma(X_2)=\sigma(X'_2)}. So the isomorphism classes of primitive double
schemes $\BS{\kx}$ such that \ \m{\BS{\kx}_{red}\simeq X} \ and with
associated line bundle $L$ are parameterized by \ \m{\P(H^1(X,T_X\ot
L))\cup\nsp}.

Another family \m{(D'_{ij})} of derivations defines (with \m{(\alpha_{ij})})
the same double scheme as \m{(D_{ij})} if and only if there exist
\m{\tau\in\C^*}, and for every $i$ a derivation \m{T_i} of \m{\ko_X(U_i)} such
that \ \m{D'_{ij}=\tau D_{ij}+T_i-\alpha_{ij}T_j}.
\end{subsub}

\end{sub}

\sepsub

\Ssect{Extension of vector bundles}{ext_hig}

Suppose that $X$ is a projective scheme. We have a canonical exact sequence of
coherent sheaves on $X$
\[0\lra L\lra\Omega_{\BS{\kx}|X}\lra\Omega_X\lra 0 \ , \]
corresponding to \ \m{\sigma\in\Ext^1_X(\Omega_X,L)}. Let $E$ be a vector
bundle on $X$. Then \m{\sigma} induces a linear map
\[\xi_{\BS{\kx}}:\Ext^1_X(E,E\ot\Omega_X)\lra\Ext^2_X(E,E\ot L) \ . \]

\sepprop

In \cite{dr10}, an element \ \m{\nabla(E)\in\Ext^2_X(E,E\ot L)} \ is
constructed, such that $E$ can be extended to a vector bundle on \m{\BS{\kx}}
if and only if \ \m{\nabla(E)=0}. From \cite{dr10}, theorem 7.1.2, we have \
\m{\nabla(E)=\xi_{\BS{\kx}}(\nabla_0(E))} (where \m{\nabla_0(E)} is the
canonical class of $E$, cf. \cite{dr10}, 3).

\sepprop

In particular, if \m{\BS{\kx}} is the trivial double scheme, then \
\m{\xi_{\BS{\kx}}=0}, hence every vector bundle on $X$ can be extended to a
vector bundle on \m{X_2}.

\sepprop

\begin{subsub}\label{prop11}{\bf Proposition: } Let $E$, $F$ be vector bundles
on $X$. Then \m{E\oplus F} can be extended to a vector bundle on \m{X_2} if and
only if $E$ and $F$ can.
\end{subsub}
\begin{proof}
We have \ \m{\nabla(E\oplus F)=\xi_{X_2}(\nabla_0(E\oplus F))}, where \Nligne
\m{\nabla_0(E\oplus F)\in\Ext^1_X(E\oplus F,(E\oplus F)\ot\Omega_X)} \ is the
canonical class of \m{E\oplus F}. We have
\[\Ext^1_X(E\oplus F,(E\oplus F)\ot\Omega_X)=\Ext^1_X(E,E\ot\Omega_X)\oplus
\Ext^1_X(E,F\ot\Omega_X)\oplus\qquad\qquad\]
\[\qquad\qquad\qquad\qquad\qquad\Ext^1_X(F,E\ot\Omega_X)\oplus
\Ext^1_X(F,F\ot\Omega_X) \ , \]
and proposition \ref{prop11} follows that from the construction of the
canonical class in \cite{dr10} we have:
\begin{enumerate}
\item[--] the components of \m{\nabla_0(E\oplus F)} in
\m{\Ext^1_X(E,F\ot\Omega_X)}
and \m{\Ext^1_X(F,E)\ot\Omega_X)} vanish;
\item[--] the components of \m{\nabla_0(E\oplus F)} in
\m{\Ext^1_X(E,E\ot\Omega_X)} is \m{\nabla_0(E)}, and its component in
\m{\Ext^1_X(F,F\ot\Omega_X)} is \m{\nabla_0(F)}.
\end{enumerate}

\end{proof}
\end{sub}

\sepsec

\section{Coherent sheaves on primitive double schemes}\label{CS_DB}

Let $X$ be a smooth connected scheme, and \m{X_2} a primitive double scheme
with support $X$. Let \ \m{L=\ki_{X,X_2}} \ be the associated line bundle on
$X$.

\sepsub

\Ssect{Canonical subsheaves and quotients -- Extensions of sheaves on $X$}{SQ}

Let $Z$ be a noetherian scheme and $\ke$ be a coherent sheaf on \m{X_2\times
Z}. Let \ \m{p_1:X_2\times Z\to X_2} \ be the first projection, and \
\m{\L=p_1^*(L)=\ki_{X\times Z,X_2\times Z}}, which is a line bundle on
\m{X\times Z}. Let \ \m{\ke_1=\L\ke}, and \ \m{\ke^{(1)}\subset\ke} \ the
subsheaf of of elements annihilated by $\L$. We have
\begin{enumerate}
\item[--]  $\ke/\ke_1=\ke_{|X\times Z}$ \ and \
\m{\ke_1\simeq(\ke/\ke^{(1)})\ot\L} .
\item[--] $\ke_1$, \m{\ke^{(1)}} are sheaves on \m{X\times Z}.
\end{enumerate}

Let \m{E\subset\ke} be a subsheaf of support contained in \m{X\times Z}, and \
\m{F=\ke/E}. Then $F$ is concentrated on \m{X\times Z} if and only \
\m{\ke_1\subset E\subset\ke^{(1)}}. In this case, the canonical morphism \
\m{\ke\ot\L\to\ke} \ induces a morphism
\[\phi_{\ke,E}:F\ot\L\lra E \]
such that \ \m{\imm(\phi_{\ke,E})=\ke_1}. We have \ \m{E=\ke^{(1)}} \ if and
only \m{\phi_{\ke,E}} is injective. We will note more simply
\[\phi_\ke \ = \ \phi_{\ke,\ke_1}:\ke_{|X}\ot\L\lra\ke_1 \ . \]

\sepprop

\begin{subsub} Extensions of sheaves on \m{X\times Z} -- \rm Let $A$, $B$ be
coherent sheaves on \m{X\times Z}. Now we define a map
\[\BS{\Phi}_{A,B}:\Ext^1_{X_2\times Z}(B,A)\lra\Hom(B\ot\L,A) \ .\]
Let \ \m{\sigma\in\Ext^1_{X_2\times Z}(B,A)}, and
\[0\lra A\lra\ke\lra B\lra 0\]
be the corresponding extension. Then \ \m{
\BS{\Phi}_{A,B}(\sigma)=\phi_{\ke,A}}.
\end{subsub}

\sepprop

\begin{subsub}\label{prop4}{\bf Proposition: } \m{\BS{\Phi}_{A,B}} is a linear
map, and \ \m{\ker(\BS{\Phi}_{A,B})=\Ext^1_{X\times Z}(B,A)}.
\end{subsub}
\begin{proof} Let \ \m{\sigma,\sigma'\in\Ext^1_{X_2\times Z}(B,A)}. Let
\[0\lra A\lra\kf\lra B\oplus B\lra 0\]
be the extension corresponding to \
\m{(\sigma,\sigma')\in\Ext^1_{X_2}(B\oplus B,A)}. Let
\[\phi_{\kf,A} \ = \ (\psi,\psi'): (B\ot\L)\oplus(B\ot\L)\lra A \ , \]
\m{0\to A\to\ke\to B\to 0} \ be the extension corresponding to \m{\sigma} and
\ \m{\alpha=\begin{pmatrix}I_B\\ 0\end{pmatrix}:B\to B\oplus B}. We have a
commutative diagram
\xmat{0\ar[r] & A\ar[r]\fleq[d] & \ke\ar[r]\ar[d] & B\ar[r]\ar[d]^\alpha & 0\\
0\ar[r] & A\ar[r] & \kf\ar[r] & B\oplus B\ar[r] & 0}
It follows from \ref{ext_sh} that \ \m{\psi=\phi_{\ke,A}}. Similarly we have \
\m{\psi'=\phi_{\ke',A}}.

Let \ \m{0\to A\to\ke''\to B\to 0} \ be the extension corresponding to
\m{\sigma+\sigma'}, and \Nligne \m{\beta=\begin{pmatrix}I_B\\
I_B\end{pmatrix}:B\to B\oplus B}. According
to \ref{ext_sh}, we have a commutative diagram
\xmat{0\ar[r] & A\ar[r]\fleq[d] & \ke''\ar[r]\ar[d] & B\ar[r]\ar[d]^\beta
& 0 \\
0\ar[r] & A\ar[r] & \kf\ar[r] & B\oplus B\ar[r] & 0}
It follows that \ \m{\BS{\Phi}_{A,B}(\sigma+\sigma')=\Psi+\psi'=
\BS{\Phi}_{A,B}(\sigma)+\BS{\Phi}_{A,B}(\sigma')}. The proof that \Nligne
\m{\BS{\Phi}_{A,B}(\lambda\sigma)=\lambda\BS{\Phi}_{A,B}(\sigma)} \ for every
\m{\lambda\in\C} is similar. The second assertion of proposition \ref{prop4} is
obvious.
\end{proof}

\end{sub}

\sepsub

\Ssect{Quasi locally free sheaves}{QF}

\begin{subsub}\label{QF2} Quasi locally free sheaves -- \rm
Let $P$ be a closed point of $X$. Let $M$ be a \m{\ko_{X_2,P}}-module of finite
type. Then $M$ is called {\em quasi free} if there exist non negative integers
\m{m_1,m_2} and an isomorphism \ \m{M\simeq(\ko_{X,P}\ot\C^{m_1})\oplus
(\ko_{X_2,P}\ot\C^{m_2})}~. The integers \m{m_1,m_2} are uniquely determined,
and \m{(m_1,m_2)} is called the {\em type} of the module $M$.
\end{subsub}

Let $\ke$ be a coherent sheaf on \m{X_2}. We say that $\ke$ is {\em quasi free
at $P$} if \m{\ke_P} is quasi free, and that $\ke$ is {\em quasi locally free}
on a nonempty open subset \m{U\subset X} if it is quasi free at every point of
$U$. In this case the types of the modules \m{\ke_x}, \m{x\in U}, are the same,
and for every \m{x\in U} there exists a neighborhood \m{V\subset U} of $x$ and
an isomorphism
\begin{equation}\label{equ20}\ke_{|V} \ \simeq \ (\ko_V\ot\C^{m_1})\oplus
(\ko_{V^{(2)}}\ot\C^{m_2}) \ .
\end{equation}
If \m{U=X} we say that $\ke$ is {\em quasi locally free}.

Let $\ke$ be a quasi locally free sheaf on \m{X_2}. Then the canonical morphism
\ \m{\ke\ot L\to\ke_1} \ induces a surjective morphism \
\m{\BS{\Psi}_\ke:\ke_{|X}\to\ke_1\ot L^*} \ of vector bundles on $X$.

The proof of the following theorem is the same as in the case \m{\dim(X)=1}
(th. 5.1.3 of \cite{dr2}).

\sepprop

\begin{subsub}\label{theo1}{\bf Theorem: } The following two assertions are
equivalent:
\begin{enumerate}
\item[(i)] The $\ko_{X_2,P}$-module $M$ is quasi free.
\item[(ii)] $M_1$ and \m{M/M_1} are free $\ko_{X,P}$-modules.
\end{enumerate}
\end{subsub}

\sepprop

It follows that

\sepprop

\begin{subsub}\label{theo2}{\bf Theorem: } The following two assertions are
equivalent:
\begin{enumerate}
\item[(i)] $\ke$ is quasi locally free.
\item[(ii)] $\ke_1$ and \m{\ke_{|X}} are locally free on $X$.
\end{enumerate}
\end{subsub}

\sepprop

It follows from theorem \ref{theo2} that for every coherent sheaf $\ke$ on
\m{X_2} there exists a nonempty open subset \m{U\subset X} such that \m{\ke} is
quasi locally free on $U$. The type of the modules \m{\ke_x}, \m{x\in U}, is
called the type of $\ke$.

\sepprop

If $\ke$ is quasi locally free, we have the exact sequence \
\m{0\to\ke_1\to\ke\to\ke_{|X}\to 0}, with \m{\ke_1} and \m{\ke_{|X}} locally
free on $X$. Conversely, if we have an exact sequence
\[0\lra\Gamma\lra\ke\lra\Lambda\lra 0 \ , \]
where $\Gamma$, $\Lambda$ are vector bundles on $X$, it does not follow that
$\ke$ is quasi locally free. For example, if $X$ is a surface and \m{C\subset
X} is a curve, the ideal sheaf \m{\ke=\ki_{C,X_2}} of $C$ in \m{X_2} in not
quasi locally free, but we have an exact sequence
\[0\lra\ke^{(1)}=L\lra\ke\lra\ke/\ke^{(1)}=\ko_X(-C)\lra 0 \ . \]

\sepprop

\begin{subsub}\label{prop2}{\bf Proposition: } Let $\Gamma$, $\Lambda$ be
vector bundles on $X$,
\[0\lra\Gamma\hfl{\iota}{}\ke\hfl{\rho}{}\Lambda\lra 0 \]
an exact sequence, where $\ke$ is a coherent sheaf on \m{X_2},
\m{\phi:\Lambda\ot L\to\Gamma} \ the associated morphism, and \m{x\in X}. Then
$\ke$ is quasi free at $x$ if and only if $\phi$ has a maximal rank at $x$.
\end{subsub}
\begin{proof}
We have \ \m{\imm(\phi)=\ke_1\subset\Gamma}. Hence $\phi$ has a maximal rank at
$x$ if and only if \m{\ke_1} is a sub-bundle of $\Gamma$ on a neighborhood of
$x$, if and only if \m{W=\Gamma/\ke_1} is locally free at $x$. We have a
commutative diagram with exact rows
and columns
\xmat{& & & 0\ar[d] \\ & 0\ar[d] & & W\ar[d] \\
0\ar[r] & \ke_1\ar[r]\ar[d] & \ke\ar[r]\fleq[d] & \ke_{|X}\ar[r]\ar[d] & 0\\
0\ar[r] & \Gamma\ar[r]\ar[d] & \ke\ar[r] & \Lambda\ar[r]\ar[d] & 0\\
& W\ar[d] & & 0\\ & 0 & & }
which shows that \m{\ke_1} and \m{\ke_{|X}} are locally free at $x$ if and only
if $W$ is locally free at $x$.
\end{proof}

\sepprop

It follows that $\ke$ is quasi locally free if and only $\phi$ has constant
rank.

\end{sub}

\sepsub

\Ssect{Some useful $\Ext$ sheaves}{DCS}

Let $\ke$ be a coherent sheaf on \m{X_2}. Then \m{\Ext^1_{X_2}(\ke,\ke)}
is the tangent space of the versal infinitesimal deformation of $X$ (cf.
\cite{si-tr}). We have the canonical Ext exact sequence
\begin{equation}\label{equ4}
0\lra H^1(X_2,\HHom(\ke,\ke))\lra\Ext^1_{X_2}(\ke,\ke)\hfl{\gamma}{}
H^0(X_2,\EExt^1_{X_2}(\ke,\ke)) \qquad\qquad\qquad\qquad\end{equation}
\[\hskip 180pt \lra H^2(X_2,\HHom(\ke,\ke)) \ . \]
Hence it is interesting to compute \m{\EExt^1_{X_2}(\ke,\ke)}. More
generally, we will prove

\sepprop

\begin{subsub}\label{theo5}{\bf Theorem: } Let $\ke$, $\kf$ be quasi locally
free sheaves on \m{X_2}. Then there is a canonical functorial isomorphism
\[\EExt^1_{X_2}(\ke,\kf) \ \simeq \ \HHom((\ke^{(1)}/\ke_1)\ot
L,\kf^{(1)}/\kf_1) \ . \]
\end{subsub}

There are integers \m{a_1,a_2\geq 0} such that $\ke$ is locally isomorphic to
\ \m{(\ko_{X_2}\ot\C^{a_2})\oplus(\ko_X\ot\C^{a_1})}.

Let $U$ be a nonempty affine open subset of $X$ such that \m{L_{|U}} is
trivial, \m{z\in\Gamma(L_{|U})} a generator, and
\[\kg \ = \ (\ko_{U^{(2)}}\ot\C^{a_2})\oplus(\ko_U\ot\C^{a_1}). \]
Suppose that there is an isomorphism \ \m{\theta:\ke_{|U^{(2)}}\to\kg}. Let \
\m{\E=\ko_{U^{(2)}}\ot\C^{a_2+a_1}} \ and \
\m{\pi:\E\to\kg} \ the canonical surjective morphism. We have
\[\ker(\pi) \ = \ L_{|U}\ot\C^{a_1} \ = \ L_{|U}\ot(\kg^{(1)}/\kg_1) . \]
It follows that we have an isomorphism
\[\eta_\theta:L_{|U}\ot(\ke_{|U}^{(1)}/\ke_{1|U})\lra\ker(\pi) \ . \]
Let \ \m{v\in\Gamma(\ke_{|U}^{(1)}/\ke_{1|U})} and \
\m{\ov{v}\in\Gamma(\ke_{|U}^{(1)})} \ over $v$. Then we can write \
\m{\theta(\ov{v})=zw_2+w_1}, \ \m{w_i\in\Gamma(\ko_U\ot\C^{a_i})}. We have then
\[\eta_\theta(z\ot v) \ = \ z\ot w_1 \ . \]
Let \ \m{\iota:L_{|U}\ot\C^{a_1}\to\E} \ be the inclusion. Then we have an
exact sequence
\begin{equation}\label{equ3}0\lra
L_{|U}\ot(\ke_{|U}^{(1)}/\ke_{1|U})\hfl{\lambda_\theta}{}
\E\hfl{\mu_\theta}{}\ke_{|U^{(2)}}\lra 0 \ , \end{equation}
with \ \m{\lambda_\theta=\iota\circ\eta_\theta} \ and \
\m{\mu_\theta=\theta^{-1}\circ\pi}.

Let $\phi$ be an automorphism of $\kg$. Then there exists an automorphism
$\ov{\phi}$ of $\E$ such that \ \m{\phi\circ\pi=\pi\circ\ov{\phi}}. Let \
\m{\theta'=\phi\circ\theta:\ke_{|U^{(2)}}\to\kg}.

\sepprop

\begin{subsub}\label{prop5}{\bf Proposition: } The following diagram
\xmat{0\ar[r] & L_{|U}\ot(\ke_{|U}^{(1)}/\ke_{1|U})\ar[r]^-{\lambda_\theta}
\fleq[d] & \E\ar[r]^-{\mu_\theta}\ar[d]^{\ov{\phi}} &
\ke_{|U^{(2)}}\ar[r]\fleq[d] & 0\\
0\ar[r] & L_{|U}\ot(\ke_{|U}^{(1)}/\ke_{1|U})\ar[r]^-{\lambda_{\theta'}} &
\E\ar[r]^-{\mu_{\theta'}} & \ke_{|U^{(2)}}\ar[r] & 0 }
is commutative.
\end{subsub}
\begin{proof} We first prove that \ \m{\ov{\phi}\circ\lambda_{\theta}=
\lambda_{\theta'}} . Let \ \m{v\in\Gamma(\ke_{|U}^{(1)}/\ke_{1|U})} \ and \
\m{\ov{v}\in\Gamma(\ke_{|U}^{(1)})} \ over $v$. Then we have \
\m{\theta(\ov{v})=zw_2+w_1}, with \m{w_i\in\Gamma(\ko_U\ot\C^{a_i})}, and
\[\eta_\theta(z\ot v) \ = \ \lambda_{\theta}(z\ot v) \ = \ z\ot w_1 \ . \]
We can write \ \m{\dsp\phi=\begin{pmatrix} \alpha & z\beta\\ \gamma & \delta
\end{pmatrix}} , with \ \m{\alpha\in\Aut(\ko_{U^{(2)}}\ot\C^{a_2})}, \m{\beta\in
\Hom(\ko_U\ot\C^{a_1},\ko_U\ot\C^{a_2})},
\Nligne\m{\gamma\in\Hom(\ko_U\ot\C^{a_2},
\ko_U\ot\C^{a_1})}, \m{\delta\in\Aut(\ko_{U}\ot\C^{a_1})}. Hence \
\m{\dsp\ov{\phi}=\begin{pmatrix} \alpha & z\beta\\ \gamma+z\gamma_1 & \delta+
z\delta_1\end{pmatrix}} , with \ \m{\gamma_1\in\Hom(\ko_U\ot\C^{a_2},
\ko_U\ot\C^{a_1})}, \m{\delta_1\in\End(\ko_U\ot\C^{a_1})}. We have
\[\ov{\phi}\circ\lambda_{\theta}(z\ot v) \ = \ z\ot\delta(w_1) \ . \]
On the other hand,
\[\theta'(\ov{v}) \ = \ \phi\circ\theta(\ov{v}) \ =
\big(z(\alpha(w_2)+\beta(w_1)), \ \delta(w_1)+z\gamma(w_2)\big) \ , \]
hence \
\[\lambda_{\theta'}(z\ot v) \ = \ z\ot(\delta(w_1)+z\gamma(w_2)) \ =
z\ot\delta(w_1) \ = \ \ov{\phi}\circ\lambda_{\theta}(z\ot v) \ , \]
\ie \m{\ov{\phi}\circ\lambda_{\theta}=
\lambda_{\theta'}} .

We have
\[\mu_{\theta'}\circ\ov{\phi} \ = {\theta'}^{-1}\circ\pi\circ\ov{\phi} \ =
\theta^{-1}\circ\phi^{-1}\circ\pi\circ\ov{\phi} \ =
\theta^{-1}\circ\phi^{-1}\circ\phi\circ\pi \ = \ \theta^{-1}\circ\pi \ = \
\mu_\theta \ . \]
\end{proof}

\sepprop

\begin{subsub}\label{Prfth5} Proof of theorem \ref{theo5} -- \rm Suppose that
we have also an isomorphism \Nligne \m{\kf_{|U^{(2)}}\simeq
(\ko_{U^{(2)}}\ot\C^{m_2})\oplus(\ko_U\ot\C^{m_1})}.
From $(\ref{equ3})$, we have an exact sequence
\[\HHom(\E,\kf_{|U^{(2)}})\hfl{\alpha}{}\HHom(L_{|U}\ot
(\ke_{|U}^{(1)}/\ke_{1|U}),\kf_{|U^{(2)}})\lra
\EExt^1_{U^{(2)}}(\ke_{|U^{(2)}},\kf_{|U^{(2)}})\lra 0 \ , \]
in which $\alpha$, given the previous isomorphisms, is the canonical morphism
\[\HHom(\ko_{U^{(2)}}\ot\C^{a_2+a_1},(\ko_{U^{(2)}}\ot
\C^{m_2})\oplus(\ko_U\ot\C^{m_1}))\lra\qquad\qquad\qquad\qquad \]
\[\qquad\qquad\qquad\qquad\HHom(L_{|U}\ot\C^{a_1},(\ko_{U^{(2)}}\ot
\C^{m_2})\oplus(\ko_U\ot\C^{m_1})) . \]
From this follows easily an isomorphism
\[\EExt^1_{U^{(2)}}(\ke_{|U^{(2)}},\kf_{|U^{(2)}}) \ \simeq \
\coker(\alpha) \ \simeq \ \HHom((\ke^{(1)}/\ke_1)_{|U}\ot
L_{|U},(\kf^{(1)}/\kf_1)_{|U}) \ . \]
The fact that these isomorphisms stick together to give theorem \ref{theo5}
follows from proposition \ref{prop5}.
\end{subsub}

\sepprop

\begin{subsub}\label{rema1}{\bf Remark:} \rm Suppose that $X$ is a
projective curve. According to \cite{dr4}, prop. 5.1, there exists a vector
bundle $\E$ on \m{X_2} and a surjective morphism \m{\E\to\ke} inducing an
isomorphism \ \m{\E_{|X}\simeq\ke_{|X}}. This is not true if \ \m{\dim(X)\geq
2}, because this would imply that every vector bundle on $X$ can be extended to
a vector bundle on \m{X_2}, which is not true in general. For example, on the
non projective double plane, the only line bundle on \m{\P_2} that can be
extended is \m{\ko_{\P_2}}.
\end{subsub}

\sepprop

\begin{subsub}\label{fam_qlfs_0} The same result for families -- \rm Let $S$ be
a scheme, and \ \m{Y=X\times S}, \m{Y_2=X_2\times S}. The notions of canonical
sequences and quasi locally free sheaves on \m{Y_2} can be defined in an
obvious way, and the proof of theorem \ref{theo5} for quasi locally free
sheaves on \m{Y_2} is similar to the case of \m{X_2}.
\end{subsub}

\sepprop

\begin{subsub} A compatibility result -- \rm Let $\ke$ be a quasi locally free sheaf on \m{X_2}. From the exact sequence \ \m{0\to\ke_1\to\ke\to\ke_{|X}\to 0} \ we deduce the map
\[\HHom(\ke,\ke_{|X})\lra\EExt^1_{X_2}(\ke,\ke_1) \ , \]
which, by theorem \ref{theo5}, is
\[\delta:\EEnd(\ke_{|X})\lra\HHom((\ke^{(1)}/\ke_1)\ot L,\ke_1)=
\ke_1\ot L^*\ot(\ke^{(1)}/\ke_1)^*  \ . \]
let \ \m{\lambda:\ke_{|X}\to\ke_1\ot L^*} \ be the canonical morphism, from which we deduce
\[H^1(\lambda\ot I_{(\ke^{(1)}/\ke_1)^*}):H^1(X,\ke_{|X}\ot(\ke^{(1)}/\ke_1)^*)
\lra H^1(X,\ke_1\ot L^*\ot(\ke^{(1)}/\ke_1)^*) \ . \]

We have \ \m{\ker(\lambda)=\ke^{(1)}/\ke_1}. Hence we get \m{\tau:(\ke_{|X})^*\to(\ke^{(1)}/\ke_1)^*}, and
\[H^1(I_{\ke_{|X}}\ot\tau):H^1(X,\EEnd(\ke_{|X}))\lra H^1(X,\ke_{|X}\ot(\ke^{(1)}/\ke_1)^*) \ . \]
The proof of the following result is similar to that of theorem \ref{theo5}, using an open cover \m{(U_i)} of $X$ such that the double schemes \m{U_i^{(2)}} are trivial, and that the \m{H^1} can be computed using \v Cech cohomology.
\end{subsub}

\sepprop

\begin{subsub}\label{prop22}{\bf Proposition: } We have \ \m{H^1(\delta)=
H^1(\lambda\ot I_{(\ke^{(1)}/\ke_1)^*})\circ H^1(I_{\ke_{|X}}\ot\tau)} .
\end{subsub}

\end{sub}

\sepsub

\Ssect{Properties of families of quasi locally free sheaves}{fam_QLF}

The following result is proved in \cite{dr12} (theorem 4.5.1):

\sepprop

\begin{subsub}\label{theo4x} {\bf Theorem: } Let $C$ be an irreducible smooth
curve. Let $\ke$ be a coherent sheaf on \m{X_2\times C}, flat on $C$. Suppose
that for every closed point \m{c\in C}, \m{\ke_c} is quasi locally free, and
that the type of \m{\ke_c} is independent of $c$. Then

{\bf 1 -- } $\ke$ is quasi locally free.

{\bf 2 -- } The sheaves \m{\ke_1} and \m{\ke_{|X\times C}} are flat on $C$, and
for every \m{c\in C}, the morphism \ \m{(\ke_1)_c\to\ke_c}
\ induced by \ \m{\ke_1\subset\ke} is injective, and
\[(\ke_1)_c \ = \ (\ke_c)_1 \ , \quad \ (\ke_{|X})_c \ =
\ke_{c|X} \ . \]

{\bf 3 -- } The sheaf \m{\ke^{(1)}} is flat on $C$, and for every \m{c\in
C},the morphism \ \m{(\ke^{(1)})_c\to\ke_c} \ induced by \
\m{\ke^{(1)}\subset\ke} is injective, and \
\m{(\ke^{(1)})_c \ = (\ke_c)^{(1)}}.
\end{subsub}

\sepprop

Let $\ke$ be a coherent sheaf on \m{X_2\times C}, flat on $C$, such that for
every closed point \m{c\in C}, \m{\ke_c} is a vector bundle on $C$. Theorem
\ref{theo4x} implies that $\ke$ is a vector bundle on \m{X\times C}.
The situation is much more complicated if we allow flat families of sheaves
parameterized by non reduced schemes. We will see in \ref{D_shX} that in
general
the deformations of $E$ parameterized by $S$ need not be sheaves on \m{X\times
S}.

\end{sub}

\sepsub

\Ssect{Deformation of sheaves on $X$}{D_shX}

Let $\ke$ be a coherent sheaf on \m{X_2}, and consider the exact sequence
$(\ref{equ4})$ in \ref{DCS}. Let \ \m{\sigma\in\Ext^1_{X_2}(\ke,\ke)}, and \
\m{0\to\ke\to\kf\to\ke\to 0} \ the corresponding extension. Let \m{x\in X} be a
closed point. Then the preceding extension at $x$ gives an exact sequence of
\m{\ko_{X_2,x}}-modules: \m{0\to\ke_x\to\kf_x\to\ke_x\to 0}. To this sequence
corresponds \
\m{\sigma_x\in\Ext^1_{\ko_{X_2,x}}(\ke_x,\ke_x)=\EExt^1_{X_2}(\ke,\ke)_x} .
We have \ \m{\gamma(\sigma)_x=\sigma_x} (cf. $(\ref{equ4})$).

\begin{subsub} Deformations of \m{\ko_{X,x}} -- \rm Let \m{t\in\ko_{X_2,x}} a
generator of the ideal of $X$, so that we have \
\m{\ko_{X_2,x}=\ko_{X,x}[t]/(t^2)}. From theorem \ref{theo5}, or from the free
resolution
\[\cdots\ko_{X_2,x}\hfl{\times t}{}\ko_{X_2,x}\hfl{\times t}{}
{}\ko_{X_2,x}\hfl{p}{}\ko_{X,x} \ , \]
we have
\[\Ext^1_{\ko_{X_2,x}}(\ko_{X,x},\ko_{X,x}) \ \simeq \
\Hom(\ko_{X,x}\ot(t),\ko_{X,x}) \ \simeq \ \ko_{X,x} \ . \]
Let \ \m{\phi\in\Ext^1_{\ko_{X_2,x}}(\ko_{X,x},\ko_{X,x})} \ and
 \begin{equation}\label{equ6}0\lra\ko_{X,x}\hfl{\iota_\phi}{}M_\phi
\hfl{\pi_\phi}{}\ko_{X,x}\lra 0\end{equation}
the corresponding extension of
\m{\ko_{X_2,x}}-modules. Then \m{M_\phi} is isomorphic to the cokernel of
\[\xymatrix@R=5pt{\theta_\phi:\ko_{X,x}\ar[r] & \ko_{X_2,x}\oplus\ko_{X,x}\\
u\fmaps[r] & (tu,\phi u)}\]
(cf. \ref{ext_sh}), \m{\iota_\phi} is induced by
\[\xymatrix@R=5pt{\ko_{X,x}\ar[r] & \ko_{X_2,x}\oplus\ko_{X,x}\\
v\fmaps[r] & (0,v)}\]
and \m{\pi_\phi} by
\[\xymatrix@R=5pt{\ko_{X_2,x}\oplus\ko_{X,x}\ar[r] & \ko_{X,x}\\
(w,u)\fmaps[r] & w_{|X}}\]
Let \ \m{{\bf Z}=\spec(\C[z]/(z^2))}. Then to $\phi$ corresponds an
infinitesimal deformation of \m{\ko_{X,x}}, \ie a sheaf on \m{X_2\times{\bf
Z}} whose restriction to \m{X_2} is \m{\ko_{X,x}}. This sheaf is \m{M_\phi}
(with an appropriate structure of \m{\ko_{X_2,x}[z]/(z^2)}-module). The
\m{\ko_{X,x}}-module structure is the same, and we need only to define the
multiplication by $z$: for every \m{u\in M_\phi}, we have \
\m{zu=\iota_\phi(\pi_\phi(u))}.
\end{subsub}

\sepprop

From this we get

\sepprop

\begin{subsub}\label{prop21}{\bf Proposition: } The sheaf \m{M_\phi} is
concentrated on \m{X\times{\bf Z}} if and only if \ \m{\phi=0}.
\end{subsub}

\sepprop

\begin{subsub}\label{coco5}{\bf Corollary: } Let $C$ be an irreducible smooth
curve. Let $\ke$ be a coherent sheaf on \m{X_2\times C}, flat on $C$. Suppose
that for every closed point \m{c\in C}, \m{\ke_c} is quasi locally free, and
that the type of \m{\ke_c} is independent of $c$. Then the image of the
Koda\"ira-Spencer morphism \ \m{TC_c\to\Ext^1_{X_2}(\ke_c,\ke_c)} \ is
contained in \ \m{H^1(X_2,\EEnd(\ke_c))}.
\end{subsub}

\end{sub}

\sepsub

\Ssect{Quasi locally free sheaves of rigid type}{QLL_R}

Let $\ke$ be a coherent sheaf on \m{X_2}, and \m{(m_1,m_2)} its type. We say
that $\ke$  is {\em of rigid type} if \m{m_1=0} or \m{m_1=1}. In this case,
for every \m{x\in X},  there exist a neighborhood $U$ of $x$ such that
\m{\ke_{|U^{(2)}}} is locally free if \m{m_1=0}, or locally isomorphic to \
\m{(\ko_{U^{(2)}}\ot\C^{m_2})\oplus\ko_U} \ if \m{m=1}.

\sepprop

\begin{subsub}\label{theo6}{\bf Theorem: } Let $Y$ be a scheme, $\ke$ a
coherent sheaf on \m{X_2\times Y}, flat on $Y$. Let \m{y_0} be a closed point
of $Y$ such that \m{\ke_{y_0}} is quasi locally free of rigid type. Then there
exists an open subset \m{U\subset Y} such that \m{y_0\in U}, and
\begin{enumerate}
\item[--] For every closed point \m{y\in U}, \m{\ke_y} is quasi locally free of
rigid type, of the same type as \m{\ke_{y_0}}.
\item[--] $\ke_{|X\times U_{red}}$ and \m{\ke_{1|X\times U_{red}}} are vector
bundles on \m{X\times U_{red}}, and for every closed point \m{y\in U}
\[(\ke_{|U_{red}\times X})_y \ = \ \ke_{y|X} \ , \quad (\ke_{1|U_{red}\times
X})_y \ =(\ke_y)_1 \ . \]
\end{enumerate}
\end{subsub}
\begin{proof}
If \m{\ke_{y_0}} is locally free, the result is well known. Suppose that
\m{\ke_{y_0}} is of type \m{(m,1)}.

By proposition \ref{prop6}, {\bf 1}, we can assume that $Y$ is integral and \
\m{\dim(Y)>0}.

Let \m{x\in X} be a closed point. Let \m{U\subset X} be an affine open subset
such that \m{x\in U}, and such that \m{U^{(2)}} is trivial, \ie \ \m{U^{(2)}
\simeq U\times{\bf Z}_2},
\m{\Gamma(\ko_{U^{(2)}})\simeq\Gamma(\ko_U)[z]/(z^2)}  (cf. \ref{PMS-1}), and
that \ \m{\ke_{y_0|U^{(2)}}\simeq(\ko_{U^{(2)}}\ot\C^{m_2})\oplus\ko_U}.

There exists a canonical surjective morphism
\[\ko_{X_2\times Y,(y_0,x)}\ot\C^{m+1}\lra\ke_{(y_0,x)} \ , \]
and a open affine subset \ \m{V\subset Y\times X} \ containing \m{(y_0,x)} such
that this morphism can be extended to a surjective morphism
\[\phi:\ko_{V^{(2)}}\ot\C^{m+1}\lra\ke_{|V^{(2)}} \ . \]
Let \ \m{\kf=\ker(\phi)}. By \cite{sga1}, Expos\'e IV, prop. 1.1, $\kf$ is
flat  on $Y$. Let \ \m{V_0=V\cap(\{y_0\}\times X)}. We have \
\m{\kf_{y_0}\simeq\ko_{V_0}}. There exists an affine open subset $W$ of $V$
containing \m{(y_0,x)} such that we have a surjective morphism \
\m{\pi:\ko_{W^{(2)}}\to\kf_{|W^{(2)}}}. We have \
\m{\Gamma(\ko_{W^{(2)}})=\Gamma(\ko_W)[z]/(z^2)} \ and
\ \m{\kf_{|W}\simeq\ko_{W^{(2)}}/\ki}, where $\ki$ is an ideal sheaf. Let
\m{\ki_0} be the ideal sheaf of \ \m{W_0=W^{(2)}\cap(X_2\times\{y_0\})}.
Since \ \m{\kf_{y_0}\simeq\ko_{V_0}}, we have
\[\ki+\ki_0 \ = \ (z)+\ki_0 \ . \]
We can write \ \m{\pi(z)=\pi(u_0)}, with \ \m{u_0\in\Gamma(\ki_0)}. We have \
\m{\pi(zu_0)=\pi(z^2)=0}. Since $\kf$ is flat on $Y$, we have \
\m{u_0\ot\pi(z)=0} \ in \ \m{\ki_0\ot_{\ko_Y}\kf_{|W^{(2)}}}. On \
\m{W^{(2)}\backslash W_0}, we have \ \m{\ki_0\simeq\ko_{W^{(2)}}}, hence
\m{u_0=0} or \m{\pi(z)=0} on \ \m{W^{(2)}\backslash W_0}. If \m{u_0=0} on \
\m{W^{(2)}\backslash W_0}, since \  \ \m{\pi(z)=\pi(u_0)} \ and \m{\pi(z)=0}
on \m{W_0}, we have \ \m{\pi(z)=0}, \ie \m{\kf_{|W^{(2)}}} is concentrated on
$W$. From proposition \ref{prop6}, {\bf 2-}, \m{\kf_{|W^{(2)}}}, as a sheaf on
$W$, is flat on $Y$. From proposition \ref{prop6}, {\bf 3-}, \m{\kf_{|W^{(2)}}}
is locally free, as a sheaf on $W$, on a neighborhood of \m{(y_0,x)}. By
shrinking $W$, we can assume that \ \m{\kf_{|W^{(2)}}\simeq\ko_W}, and that the
inclusion
\xmat{\ker(\phi)_{|W}=\ko_W\flinc[rr] & & \ko_{|W^{(2)}}\ot\C^{m+1}}
induces an injective morphism of vector bundles on $W$
\[\ko_W\lra\big[\ko_{|W^{(2)}}\ot\C^{m+1}\big]^{(1)}\simeq\ko_W\ot\C^{m+1} \ .
\]
It follows easily that \
\m{\ke_{|W^{(2)}}\simeq(\ko_{W^{(2)}}\ot\C^{m+1})\oplus\ko_W} .

There is a finite number of open subsets \m{W_1,\ldots,W_k} of \ \m{X\times
Y}, such that \Nligne \m{X\times\{y_0\}\subset\bigcup_{1\leq i\leq k}W_i}, and
for \m{1\leq i\leq k},
\m{\ke_{|W_i^{(2)}}\simeq(\ko_{W_i^{(2)}}\ot\C^{m+1})\oplus\ko_{W_i}} . Let
\Nligne \m{p_Y:X\times Y\to Y} \ be the projection, and
\[Z \ = \ p_Y\big((X\times Y)\backslash(\cup_{1\leq i\leq k}W_i)\big) \ . \]
It is a closed subset of $Y$, and for every \ \m{y\in Y\backslash Z}, \m{\ke_y}
is quasi locally free of type \m{(1,m)}.
\end{proof}

\end{sub}

\sepsub

\Ssect{Isomorphisms of quasi locally free sheaves}{Isom}

Let $\E$, $E$ be simple vector bundles on $X$. Let $\ke$, $\ke'$ be quasi
locally free sheaves on \m{X_2}, with isomorphisms
\[\theta_X:\ke_{|X}\to\E \ , \quad \theta'_X:\ke'_{|X}\to\E \ , \quad
\theta_1:\ke_1\to E\ot L \ , \quad \theta'_1:\ke'_1\to E\ot L \ , \]
from which we deduce
\begin{enumerate}
\item[--] the exact sequences
\[0\lra E\ot L\lra\ke\lra\E\lra 0 \ , \qquad 0\lra E\ot L\lra\ke'\lra\E\lra 0 \]
associated to \ \m{\sigma,\sigma'\in\Ext^1_{X_2}(\E,E\ot L)} \ respectively,
\item[--] the induced morphisms \ \m{\Phi:\E\to E}, \m{\Phi':\E\to E} .
\end{enumerate}

\sepprop

\begin{subsub}\label{prop12}{\bf Proposition: } {\bf 1 -- } We have \
\m{\ke\simeq\ke'} \ if and only if there exists \m{\lambda\in\C^*} such that \
\m{\sigma'=\lambda\sigma} . In this case we have \
\m{\dsp\Phi'=\frac{1}{\lambda}\Phi} .

{\bf 2 -- } Suppose that the map
\[\xymatrix@R=5pt{H^0(X,L)\ar[r] & \Hom(\E,E\ot L)\\
\lambda\fmaps[r] & (I_E\ot\lambda)\circ\phi}\]
is surjective. Then $\ke$ is simple (\ie \ \m{\End(\ke)=H^0(X_2,\ko_{X_2})}).
\end{subsub}
\begin{proof}
Let \ \m{\alpha:\ke\to\ke'} \ be an isomorphism. Then $\alpha$ induces
isomorphisms \Nligne \m{\alpha_1:\ke_1\ot L\to\ke'_1\ot L},
\m{\alpha_{|X}:\ke_{|X}\to\ke'_{|X}}, and we have a commutative diagram
\xmat{0\ar[r] & \ke_1\ar[r]\ar[d]^{\alpha_1} & \ke\ar[r]\ar[d]^\alpha &
\ke_{|X}\ar[r]\ar[d]^{\alpha_{|X}} & 0\\
 0\ar[r] & \ke'_1\ar[r] & \ke'\ar[r] & \ke'_{|X}\ar[r] & 0}
whence a commutative diagram
\xmat{0\ar[r] & E\ot L\ar[r]\ar[d]^{\beta_1} & \ke\ar[r]\ar[r]\ar[d]^\alpha &
\E\ar[r]\ar[d]^{\beta_X} & 0\\
0\ar[r] & E\ot L\ar[r] & \ke\ar[r]\ar[r] & \E\ar[r] & 0}
with \ \m{\beta_1=\theta'_1\circ\alpha_1\circ\theta_1^{-1}},
\m{\beta_X=\theta'_X\circ\alpha_{|X}\circ\theta_X^{-1}}. Then {\bf 1} follows
immediately from the fact that \m{\beta_1} and \m{\beta_X} are homotheties,
since $E$ and $\E$ are simple. The proof of {\bf 2} is similar.
\end{proof}

\sepprop

It follows that the isomorphism classes of quasi locally free sheaves \m{\ke'}
on \m{X_2} such that \ \m{\BS{\Psi}_{\ke'}\simeq\BS{\Psi_\ke}} \ form an affine
space with associated vector space \ \m{\Ext^1_X(\E,E\ot L)}.

\end{sub}

\sepsec

\section{Construction of quasi locally free sheaves on primitive double
schemes}\label{QLF}

Let $X$ be an irreducible smooth variety, \m{X_2} a double primitive
scheme such that \ \m{(X_2)_{red}=X}, and with associated line bundle \
\m{\ki_{X,X_2}=L}.

If \m{U\subset X} is an open subset, let \m{U^{(2)}} denote
the corresponding open subscheme of \m{X_2}. We say that {\em \m{X_2} is
trivial on $U$} if there is an isomorphism \ \m{U^{(2)}\simeq U\times{\bf Z}_2}
(with \ \m{{\bf Z}_2=\spec(\C[t]/(t^2))}) inducing the identity on $U$.

Let \ \m{\Phi:\E\to E} , \m{\Phi':\E'\lra E'} \ be surjective morphisms of
vector bundles on $X$. We say that $\Phi$ and \m{\Phi'} are {\em isomorphic}
(or \ \m{\Phi\simeq\Phi'}) if there exist isomorphisms \ \m{\alpha:\E\to\E'},
\m{\beta:E\to E'} \ such that \ \m{\Phi'\circ\alpha=\beta\circ\Phi}.

\sepsub

\Ssect{The obstruction to the construction of quasi locally free sheaves}{oQLL}

Let $F$ (resp. $\kf$) be vector bundles on $X$.

\begin{subsub}\label{ext_ss} The canonical sequence of extensions -- \rm We
have a functorial exact sequence
\begin{equation}\label{equ5}0\lra\Ext^1_X(F,\kf\ot L))\lra\Ext^1_{X_2}(F,\kf\ot
L)\hfl{\gamma_{F,\kf}}{}\Hom(F,\kf) \qquad\qquad
\end{equation}
\[\qquad\qquad\qquad \hfl{\nabla_{F,\kf}}{}\Ext^2_X(F,\kf\ot L)\]
This can be proved easily with \v Cech cohomology, using \ref{cech} and
\ref{PMS-1}. This is also a consequence of a more general result, the {\em Ext
spectral sequence} (cf. \cite{go}, 7.3 and lemma \ref{lem2}). In this case we
use the fact that we have \ \m{\EExt^1_{X_2}(F,\kf\ot L)\simeq\HHom(F,\kf)} \
by theorem \ref{theo5}.

We call this exact sequence the {\em canonical exact sequence of extensions}.
\end{subsub}

\sepprop

Let $\ke$ be a coherent sheaf on \m{X_2}. Let
\[\BS{\Psi}_\ke=\phi_\ke\ot I_{L^*}:\ke_{|X}\lra\ke_1\ot L^*\]
be the canonical morphism (cf. \ref{SQ}).

We will prove that there exists a quasi locally free sheaf $\ke$ on \m{X_2}
such that \ \m{\Phi\simeq\BS{\Psi}_\ke} \ if and only if \
\m{\nabla_{\E,E}(\Phi)=0}.

In the case of vector bundles, i.e \m{\E=E} and \m{\phi=I_E}, there is already
an invariant in \m{\Ext^2_X(E,E\ot L)}, \m{\nabla(E)}, that vanishes if and
only $E$ can be extended to a vector bundle on \m{X_2}, cf. \ref{ext_hig}.
Using its construction in \cite{dr10} (which is a particular case of what is
done in \ref{S_G}) and the construction of $(\ref{equ5})$ by \v Cech
cohomology, it is easy to see that these two invariants are the same, \ie \
\m{\nabla(E)=\nabla_{\E,\E}(I_\E)}.

Let
\[\eta_\Phi:\Ext^2_X(\E,\E\ot L)\lra\Ext^2_X(\E,E\ot L) \ , \qquad
\mu_\Phi:\Ext^2_X(E,E\ot L)\lra\Ext^2_X(\E,E\ot L)\]
be the maps induced by $\Phi$. We will prove

\sepprop

\begin{subsub}\label{theo3}{\bf Theorem:} We have \
\m{\nabla_{\E,E}(\Phi)=\eta_\Phi(\nabla(\E))=\mu_\Phi(\nabla(E))} .
\end{subsub}

\sepprop

There is another way to prove theorem \ref{theo3}. It uses sheaves of algebraic
groups associated to $\Phi$. We will follow the first way here and the second
in \ref{S_G}.

\end{sub}

\sepsub

\Ssect{Proves using the canonical exact sequence of extensions}{pr_ext}

Let $F$ (resp. $\kf$) be a vector bundle (resp. a coherent sheaf) on $X$.

\sepprop

\begin{subsub}\label{lem2}{\bf Lemma: } We have a canonical isomorphism \
\m{\EExt^1_{X_2}(F,\kf\ot L)\simeq\HHom(F,\kf)}.
\end{subsub}
\begin{proof}
If $\kf$ is quasi locally free, this follows from theorem \ref{theo5}.
Let \m{U\subset X} be an affine open subset of $X$. Since \m{U^{(2)}} is
trivial, \m{F_{|U}} and \m{L_{|U}} can be extended to vector bundles $\F$ and
$\L$ on \m{U^{(2)}}. We have a canonical morphism \
\m{\lambda:\L\to\ko_{|U^{(2)}}}
\xmat{\L\flon[r] & \L_{|U}=L=\ki_{U,U^{(2)}}\flinc[r] & \ko_{U^{(2)}} \ . }
We have a locally free resolution of \m{F_{|U}} on \m{U^{(2)}}
\xmat{\cdots\F\ot\L^2\ar[rr]^-{I_\F\ot \lambda\ot I_\L} & &
\F\ot\L\ar[rr]^-{I_\F\ot\lambda}
 & & \F\ar[r] & F_{|U} \ , }
from which we obtain an isomorphism \
\m{\EExt^1_{U^{(2)}}(F_{|U},\kf_{|U}\ot L)\simeq\HHom(F_{|U},\kf_{|U})} .
It is easy to see that these isomorphisms can be glued together to define that
of lemma \ref{lem2}.
\end{proof}

\sepprop

\begin{subsub}\label{lem3}{\bf Lemma: } Let \
\m{\sigma\in\Ext^1_{X_2}(F,\kf\ot L)}, corresponding to the exact
sequence
\[0\lra\kf\ot L\hfl{\alpha}{}\ke\hfl{\beta}{}F\lra 0\]
of sheaves on \m{X_2}, and \ \m{p:F\to\kf} \ the associated morphism (cf.
\ref{SQ}). Then we have \ \m{\gamma_{F,\kf}(\sigma)=p}.
\end{subsub}
\begin{proof}
Let \ \m{x\in X}. Then \
\m{\gamma_{F,\kf}(\sigma)(x)\in\Ext^1_{\ko_{X_2,x}}(F_x,(\kf\ot L)_x)} \
corresponds to the exact sequence
\[0\lra (\kf\ot L)_x\hfl{\alpha_x}{}\ke_x\hfl{\beta_x}{}F_x\lra 0\]
of \m{\ko_{X_2,x}}-modules, and the fact that \
\m{\gamma_{F,\kf}(\sigma)(x)=p_x} \ follows easily from the proof of lemma
\ref{lem2} and \ref{ext_sh}.
\end{proof}

\sepprop

\begin{subsub}\label{coro1}{\bf Corollary: } Let \ \m{p:F\to\kf} \ be a
morphism. Then there exists an exact sequence of sheaves on \m{X_2}
\[0\lra\kf\ot L\lra\ke\lra F\lra 0\]
whose associated morphism \ \m{F\to\kf} \ is $p$ if and only if \
\m{\nabla_{F,\kf}(p)=0}.
\end{subsub}

\sepprop

\begin{subsub}\label{coro2}{\bf Corollary: }
Let \ \m{\Phi:\E\to E} be a surjective morphism of vector bundles on $X$. Then
there exists a quasi locally free sheaf $\ke$ on \m{X_2} such that \
\m{\Phi\simeq\BS{\Psi}_\ke} if and only if \ \m{\nabla_{\E,E}(\Phi)=0}.
\end{subsub}

\sepprop

Now we prove theorem \ref{theo3}. From the functoriality of the Ext spectral
sequence, we have a commutative diagram
\xmat{\End(E)\ar[r]^-{\nabla_{E,E}}\ar[d]^\psi & \Ext^2_X(E,E\ot
L)\ar[d]^{\mu_\Phi}\\
\Hom(\E,E)\ar[r]^-{\nabla_{\E,E}} & \Ext^2_X(\E,E\ot L)}
(where $\psi$ is induced by $\Phi$). It follows that
\m{\nabla_{\E,E}(\Phi)=\mu_\Phi(\nabla_{E,E}(I_E))}. The proof of the other
equality is similar.

\sepprop

The following result is an immediate consequence of theorem \ref{theo3}, but we
will give a direct proof.

\sepprop

\begin{subsub}\label{coro3}{\bf Corollary: } Let \ \m{\Phi:\E\to E} be a
surjective morphism of vector bundles on $X$. If $E$ or $\E$ can be extended
to a vector bundle on $X_2$, then there exists a quasi locally free sheaf $\ke$
on \m{X_2} such that \ \m{\Phi\simeq\BS{\Psi}_\ke}
\end{subsub}
\begin{proof}
Suppose that $E$ can be extended to a vector bundle $\BS{E}$ on \m{X_2}. We
have a canonical exact sequence
\[0\lra E\ot L\lra\BS{E}\lra E=\BS{E}_{|X}\lra 0 \]
Let \ \m{\sigma\in\Ext^1_{X_2}(E,E\ot L)} \ be the corresponding
element, and
\[\rho:\Ext^1_{X_2}(E,E\ot L)\lra\Ext^1_{X_2}(\E,E\ot L)\]
the map induced by $\Phi$. We have then a commutative diagram with exact rows
\xmat{0\ar[r] & E\ot L\ar[r]\fleq[d] & \ke\ar[r]\ar[d] & \E\ar[r]\flon[d]^\Phi
& 0 \\
0\ar[r] & E\ot L\ar[r] & \BS{E}\ar[r] & E\ar[r] & 0 }
where the top exact sequence is associated to \m{\rho(\sigma)} (cf.
\cite{dr1b}, 4.3). Locally, on an open subset \m{U\subset X} such that \
\m{\BS{E}_{|U^{(2)}}\simeq\ko_{U^{(2)}}\ot\C^{m+n}}, and \m{L_{|U}} is
generated by $z$, it is the canonical diagram
\xmat{0\ar[r] & z\ko_{U^{(2)}}\ot\C^m\ar[r]\fleq[d] &
(\ko_{U^{(2)}}\ot\C^m)\oplus(\ko_U\ot\C^n)\ar[r]\ar[d]^{(I,0)} &
\ko_U\ot\C^{m+n}\ar[r]\flon[d]^p
& 0 \\
0\ar[r] & z\ko_{U^{(2)}}\ot\C^m\ar[r] & \ko_{U^{(2)}}\ot\C^m\ar[r] &
\ko_U\ot\C^m\ar[r] & 0 }
Hence $\ke$ is locally isomorphic to \
\m{(\ko_{X_2}\ot\C^m)\oplus(\ko_X\ot\C^n)}, and the morphism \ \m{\E\to E} \
associated to \m{0\to E\ot L\to\ke\to\E\to 0} \ is \m{\Phi}.

The case where $\E$ can be extended to a vector bundle on $X_2$ is similar.
\end{proof}

\end{sub}

\sepsub

\Ssect{Examples -- General results}{Ex0}

 Let \ \m{\Phi:\E\to E} be a surjective morphism of vector bundles on $X$.

\sepprop

\begin{subsub}\label{prop4b}{\bf Proposition: } Suppose that $\E$ is simple,
\m{\rk(E)<\rk(\E)}, and \ \m{h^0(X,\omega_X^*\ot L)>0}. Then \
\m{\Ext^2_X(\E,E\ot L)=\nsp}.
\end{subsub}
\begin{proof}
By Serre duality, we have an isomorphism
\[\Ext^2_X(\E,E\ot L) \ \simeq \ \Hom(E,\E\ot\omega_X\ot L^*)^* \ . \]
Suppose that \ \m{\Ext^2_X(\E,E\ot L)\not=\nsp}. Then there exists a nonzero
morphism \Nligne \m{\alpha:E\to\E\ot\omega_X\ot L^*}. There exists also a
nonzero morphism \ \m{\beta:\omega_X\ot L^*\to\ko_X}. Hence \
\m{\psi=(I_\E\ot\beta)\circ\alpha} \ is a nonzero morphism \ \m{E\to\E}, and \
\m{\phi\circ\Phi:\E\to\E} \ is a nonzero morphism of rank \m{<\rk(\E)}. This is
impossible since $\E$ is simple.
\end{proof}

\sepprop

This implies that if the hypotheses of this proposition are satisfied, then
there exists a quasi locally free sheaf $\ke$ on \m{X_2} such that  \
\m{\Phi\simeq\BS{\Psi}_\ke}.

\sepprop

\begin{subsub}\label{lem5}{\bf Lemma: } Let \m{E'=\ker(\Phi)}. Suppose that \
\m{\Hom(E',E)=\nsp}, \m{\Hom(E,E')=\nsp}, and that the extension
\[0\lra E'\lra\E\hfl{\Phi}{}E\lra 0\]
is not trivial . Then $\E$ is simple.
\end{subsub}
\begin{proof} Let \ \m{\alpha\in\End(\E)}. Then since \ \m{\Hom(E',E)=\nsp},
\m{(\Phi\circ\alpha)_{|E'}=0}, and \ \m{\alpha(E')\subset E'}. Let \
\m{\alpha'=\alpha_{|E'}:E'\to E'}. The fact that \m{E'} is simple implies that
there exists \m{\lambda\in\C} such that \ \m{\alpha'=\lambda I_{E'}}. Let \
\m{\eta=\alpha-\lambda I_\E}. Then \ \m{\eta(E')=\nsp}, and $\eta$ induces a
morphism \ \m{\ov{\eta}:E\to\E} \ such that \ \m{\eta=\ov{\eta}\circ\Phi}.
Since $E$ is simple, there exists \ \m{\mu\in\C} \ such that \
\m{\Phi\circ\ov{\eta}=\mu I_{E}}. Suppose that \ \m{\mu\not=0}. Then
\m{\dsp\frac{1}{\mu}\ov{\eta}} is a section of $\Phi$, which is impossible
because the extension is not trivial. Hence \ \m{\mu=0}, and \
\m{\ov{\eta}(E)\subset E'}. But \ \m{\Hom(E,E')=\nsp}, hence \ \m{\ov{\eta}=0},
\m{\eta=0} \ and \ \m{\alpha=\lambda I_\E}.
\end{proof}

\sepprop

{\bf Remark: }Suppose that \ \m{\Ext^2_X(\E,E\ot L)=\nsp}. Let \ \m{\Phi:\E\to
E} be a surjective morphism, and  \ \m{D\in\Pic(X)}. Then we have also \
\m{\Ext^2_X(\E\ot D,E\ot D\ot L)=\nsp}. Hence there exist exact sequences on
\m{X_2}
\[0\lra (E\ot D)\ot L\lra\kf\lra\E\ot D\lra 0 \ , \]
whose associated morphism \ \m{\E\ot D\to E\ot D} \ is \m{\Phi\ot I_D}, and
with $\kf$ quasi locally free.

If $D$ can be extended to a line bundle $\BS{\kd}$ on \m{X_2}, these exact
sequences are exactly those corresponding to \m{D=\ko_X}, twisted by
$\BS{\kd}$. But they exist even if $D$ cannot be extended
to a line bundle on \m{X_2}.

\sepprop

\begin{subsub}\label{DirS}{\bf Proposition: } Let $E$, $F$ be vector bundles
on $X$ and \ \m{\E=E\oplus F}. Let \ \m{\Phi:\E\to E} \ be the projection. Then
there exists a quasi locally free sheaf $\ke$ on \m{X_2} such that \
\m{\Phi\simeq\Phi_\ke} \ if and only if $E$ can be extended to a vector bundle
on \m{X_2}.
\end{subsub}
\begin{proof} We have \ \m{\Ext^2_X(\E,E\ot L)\simeq\Ext^2_X(E,E\ot L)\oplus
\Ext^2_X(F,E\ot L)}, and
\[\mu_\Phi:\Ext^2_X(E,E\ot L)\lra\Ext^2_X(\E,E\ot L)\]
is the projection on the first factor. The result follows immediately from the
fact that \ \m{\nabla_{\E,E}(\Phi)=\mu_\Phi(\nabla(E))}, by theorem
\ref{theo3}.
\end{proof}

\end{sub}

\sepsub

\Ssect{Families of quasi locally free sheaves}{fam_qlf}

Let $S$, $T$ be smooth irreducible varieties and $\E$ (resp. $\F$) a vector
bundle on \m{X\times S} (resp. \m{X\times T}). Suppose that, for every closed
points \m{s\in S}, \m{t\in T},
\begin{enumerate}
\item[(i)] the dimensions of the spaces \m{\Hom(\E_s\ot L,\F_t)},
\m{\Ext^1_X(\E_s,\F_t)}, \m{\Ext^2_X(\E_s,\F_t)}, do not depend on \m{s,t},
\item[(ii)] there are surjective morphism \ \m{\E_s\to\F_t\ot L^*}, and for
every such morphism $\Phi$, we have \ \m{\nabla_{\E_s,\F_t\ot L^*}(\Phi)=0}.
\end{enumerate}

Let \m{p_S:X\times S\times T\to X\times S}, \m{p_T:X\times S\times T\to X\times
T}, \m{p_{S\times T}:X_2\times S\times T\to S\times T} be the projections.

\sepprop

\begin{subsub}\label{theo1b}{\bf Theorem:} \m{{\bf Ext}^1_{p_{S\times
T}}(p_S^*(\E), p_T^*(\F))} commutes with base change.
\end{subsub}
\begin{proof}
(cf. \ref{rel_e}). Here of course, \m{p_S^*(\E)} and \m{p_T^*(\F)} are seen as
sheaves on the primitive double scheme \m{X_2\times S\times T} concentrated on
\m{X\times S\times T}.

We will use the canonical exact sequence of extensions (cf. \ref{ext_ss})

Let \m{q_{S\times T}:X\times S\times T\to S\times T}, \m{p_X:X\times S\times
T\to X} be the projections. Then \Nligne \m{{\bf Ext}^1_{q_{S\times
T}}(p_S^*(\E),p_T^*(\F))} (resp. \m{{\bf Ext}^0_{q_{S\times T}}(p_S^*(\E)\ot
p_X^*(L),p_T^*(\F))}, \m{{\bf Ext}^2_{q_{S\times
T}}(p_S^*(\E),p_T^*(\F))}) is locally free on \m{S\times T}, commutes with base
change, and
\[{\bf Ext}^1_{q_{S\times T}}(p_S^*(\E),p_T^*(\F)) \ \simeq \
R^1{q_{S\times T}}_*(p_S^*(\E)^*\ot p_T^*(\F)) \ , \]
\[{\bf Ext}^0_{q_{S\times T}}(p_S^*(\E)\ot p_X^*(L),p_T^*(\F)) \ = \
\HHom(p_S^*(\E)\ot p_X^*(L),p_T^*(\F)) \ , \]
\[{\bf Ext}^2_{q_{S\times T}}(p_S^*(\E),p_T^*(\F)) \ \simeq \
R^2{q_{S\times T}}_*(p_S^*(\E)^*\ot p_T^*(\F)) \ . \]
Let \ \m{U\subset S\times T} be an open affine subset and \m{s\in S}, \m{t\in
T} be closed points. Let
\[\nabla_U=\nabla_{p_S^*(\E)_{|U},(p_T^*(\F)\ot p^*(L^*))_{|U}} \ . \]
By functoriality of the canonical
exact sequence of extensions we have a commutative diagram
\xmat{\Hom(p_S^*(\E)_{|U},[p_T^*(\F)\ot p_X^*(L^*)]_{|U})
\ar[r]^-{\nabla_U}\ar[d]^{C_{s,t}} &
\Ext^2_{X\times S\times T}(p_S^*(\E)_{|U},p_T^*(\F)_{|U})\ar[d]^{D_{s,t}}\\
\Hom(\E_s,\F_t\ot L^*)\ar[r]^-{\nabla_{\E_s,\F_s\ot L^*}=0} &
\Ext^2_X(\E_s,\F_t)}
where \m{C_{s,t}}, \m{D_{s,t}} are the canonical maps. Since $U$ is affine, we
have
\[\Ext^2(p_S^*(\E)_{|U},p_T^*(\F)_{|U}) \ = \ H^0(U,R^2{q_{S\times
T}}_*(p_S^*(\E)^*\ot p_T^*(\F)))
) \ , \]
and by the preceding diagram, if \ \m{\alpha\in
\Hom(p_S^*(\E)_{|U},[p_T^*(\F)\ot p_X^*(L^*)]_{|U})}, then \
\m{\Delta_U(\alpha)(s,t)=0}. Hence \ \m{\Delta_U(\alpha)=0}, and \m{\Delta_U=0}.

Let \ \m{\gamma_U=\gamma_{p_S^*(\E)_{|U},(p_T^*(\F)\ot p^*(L^*))_{|U}}} . We
have a commutative diagram
\xmat{0\ar[r] & \Ext^1_{X\times S\times T}(p_S^*(\E)_{|U},p_T^*(\F)_{|U})
\ar[r]\ar[d]^-{A_{s,t}} & \Ext^1_{X_2\times S\times
T}(p_S^*(\E)_{|U},p_T^*(\F)_{|U})\ar[r]^-{\gamma_U}\ar[d]^{B_{s,t}} & \cdots\\
0\ar[r] & \Ext^1_X(\E_s,\F_s)\ar[r] & \Ext^1_{X_2}(\E_s,\F_s)
\ar[r]^-{\gamma_{\E_s,\F_s\ot L^*}} & \cdots
}
\xmat{\qquad\qquad\qquad\cdots
\ar[r]^-{\gamma_U} & \Hom(p_S^*(\E)_{|U},[p_T^*(\F)\ot p_X^*(L^*)]_{|U})
\ar[d]^{C_{s,t}}\\
\qquad\qquad\qquad\cdots\ar[r]^-{\gamma_{\E_s,\F_s\ot L^*}} &
\Hom(\E_s\ot L,\F_t)}
Since \m{\nabla_U} and \m{\nabla_{\gamma_{\E_s,\F_s\ot L^*}}}
vanish, the maps \m{\gamma_U} and \m{\gamma_{\E_s,\F_s\ot L^*}} are surjective,
as well as \m{A_{s,t}} and \m{C_{s,t}}. It follows easily that \m{B_{s,t}} is
surjective. This proves theorem \ref{theo1b}.
\end{proof}

\sepprop

\begin{subsub}\label{mod-sp2} Moduli spaces -- \rm (cf. \ref{fin-mod}).
Suppose that moreover, for every \ \m{(s,t)\in S\times T}
\begin{enumerate}
\item[(iii)] $\E_s$ and \m{\F_t} are simple,
\item[(iv)] $\Hom(\E_s,\F_t)=\nsp$ .
\end{enumerate}
Let \ \m{\pmb{\P}=\P\big({\bf Ext}^1_{p_{S\times T}}(p_S^*(\E),
p_T^*(\F))\big)} \ the projective bundle over \m{S\times T}, and \Nligne
\m{q_S:X_2\times\pmb{\P}\to X\times S}, \m{q_T: X_2\times\pmb{\P}\to X\times
T} \ \m{q_{\pmb{\P}}:X_2\times\pmb{\P}\to\pmb{\P}} \ the projections. According
to \cite{la}, corollary 4.5, there exists a {\em universal extension} on \
\m{X_2\times\pmb{\P}}
\[0\lra q_T^*(\F)\ot q_{\pmb{\P}}^*(\ko_{\pmb{\P}}(1))\lra\kh\lra q_S^*(\E)\lra
0 \ . \]
Let \ \m{(s,t)\in S\times T}, and \ \m{\lambda\in\pmb{\P}_{(s,t)}}, which is a
line in \m{\Ext^1_X(\E_s,\F_t)}. Then the restriction of this exact sequence to
\ \m{X_2=X_2\times\{(s,t)\}} is the extension
\begin{equation}\label{equ1}0\lra\F_t\ot\lambda^*\lra\kh_{(s,t)}\lra\E_s\lra 0
\end{equation}
associated to the inclusion in \ \m{\Hom(\lambda,\Ext^1_X(\E_s,\F_t))=
\Ext^1_X(\E_s,\F_t\ot\lambda^*)}. Hence \m{\pmb{\P}} parameterizes all the non
split extensions of bundles \m{\E_s} by bundles \m{\F_t}, modulo the action of
\m{\C^*}.

Let \m{\pmb{\P}'\subset\pmb{\P}} be the open subset corresponding to extensions
$(\ref{equ1})$ such that the corresponding morphism \ \m{\E_s\ot
L\to\F_s\ot\lambda^*} \ is surjective. Then by proposition \ref{prop12},
\m{\pmb{\P}'_{(s,t)}} parameterizes all the quasi locally free sheaves $\ke$ on
\m{X_2} such that \ \m{\ke_{|X}\simeq\E_s} \ and \ \m{\ke_1\simeq\F_t}.
\end{subsub}

Suppose that
\begin{enumerate}
\item[--] For every \m{(s,t)\in S\times T} we have \
\m{\rk(\E_s)=\rk(\F_t)+1}.
Hence every quasi locally free sheaf $\ke$ on \m{X_2} such that \
\m{\ke_{|X}\simeq\E_s} \ and \ \m{\ke_1\simeq\F_t} \ is of rigid type (cf.
\ref{QLL_R}).
\item[--] The set $\chi$ of isomorphism classes of quasi locally free sheaves
$\ke$ on \m{X_2} such that there exists \m{(s,t)\in S\times T} such that \
\m{\ke_{|X}\simeq\E_s} \ and \ \m{\ke_1\simeq\F_t} is open.
\end{enumerate}
Then it is easy to see, using theorem \ref{theo6}, that $\pmb{\P}'$, equipped
with \m{\kh_{|X_2\times\pmb{\P}'}}, is a fine moduli space for family of
sheaves in $\chi$ parameterized by integral schemes.

\end{sub}

\sepsec

\section{Sheaves of algebraic groups associated to quasi locally free
sheaves}\label{S_G}

Let $X$ be an irreducible smooth projective variety, \m{X_2} a double primitive
scheme such that \ \m{(X_2)_{red}=X}, and with associated line bundle \
\m{\ki_{X,X_2}=L}.

\sepsub

\Ssect{Construction of quasi locally free sheaves}{const_QLF}

If \m{U\subset X} is an open subset, let \m{U^{(2)}} denote
the corresponding open subscheme of \m{X_2}. We say that {\em \m{X_2} is
trivial on $U$} if there is an isomorphism \ \m{U^{(2)}\simeq U\times{\bf Z}_2}
\ inducing the identity on $U$.

\mycomment{Let
\[\Sigma: \quad 0\lra F\lra\E\lra E\lra 0 \ , \quad
\Sigma': \quad 0\lra F'\lra\E'\lra E'\lra 0\]
be exact sequences of vector bundles on $X$. We say that $\Sigma$ and
\m{\Sigma'} are {\em isomorphic} (or \ \m{\Sigma\simeq\Sigma'}) if there exists
a commutative diagram
\xmat{0\ar[r] & E\ar[r]\ar[d] & \E\ar[r]\ar[d] & F\ar[r]\ar[d] & 0\\
0\ar[r] & E'\ar[r] & \E'\ar[r] & F'\ar[r] & 0}
where the vertical arrows are isomorphisms.}

Let $M$, $N$ be complex vector spaces, with \ \m{\dim(M)=m}, \m{\dim(N)=n}. Let
$\ke$ be a quasi locally free sheaf on \m{X_2}, locally isomorphic to \
\m{(\ko_{X_2}\ot M)\oplus(\ko_X\ot N)}. Then \m{\ke_1} and \
\m{\ke/\ke_1=\ke_{|X}} are vector bundles on $X$, of rank $m$, \m{m+n}
respectively. The canonical morphism \ \m{\ke\ot\ki_X\to\ke} \ induces a
surjective morphism \ \m{\ke_{|X}\to\ke_1\ot L^*} \ of vector bundles on $X$.

Let $\kg_{m,n}$ the sheaf of algebraic groups on $X$, such that for every open
subset \m{U\subset X}, \m{\kg_{m,n}(U)} is the group of automorphisms of \
\m{(\ko_{U^{(2)}}\ot M)\oplus(\ko_U\ot N)}. Such automorphisms are matrices
\m{\begin{pmatrix} A & B \\ C & E\end{pmatrix}}, where
\begin{enumerate}
\item[] $A\in\Aut(\ko_{U^{(2)}}\ot M)$ ,
\item[] $B\in\Hom(\ko_U\ot
N,\ko_{U^{(2)}}\ot M)=\Hom(\ko_U\ot N,L_{|U}\ot M)$ ,
\item[] $C\in\Hom(\ko_{U^{(2)}}\ot M,\ko_U\ot N)=\Hom(\ko_U\ot M,\ko_U\ot N)
=\ko_X(U)\ot\kl(M,N)$ ,
\item[] $E\in\Aut(\ko_U\ot N)$ .
\end{enumerate}

Suppose that \m{X_2} and $L$ are trivial on $U$, so there exists a generator
\m{z\in\ki_X(U)}. Then we can write
\[A \ = \ \A_0+z\A_1 \ , \qquad B \ = \ z\B \ , \]
with \ \m{\A_0\in\Aut(\ko_X(U)\ot M)}, \m{\A_1\in\ko_X(U)\ot\End(M)},
\m{\B\in\ko_X(U)\ot \kl(N,M)}. We have
\[\begin{pmatrix} A & B \\ C & E\end{pmatrix}^{-1} \ = \ \begin{pmatrix}
\A_0^{-1}+
z\A_0^{-1}\B E^{-1}C\A_0^{-1} & -z\A_0^{-1}\B E^{-1} \\ -E^{-1}C\A_0^{-1} &
E^{-1}\end{pmatrix} \ . \]

The cohomology set \m{H^1(X,\kg_{m,n})} (cf. \ref{coh_gr}) is naturally
identified with the set of isomorphism classes of quasi locally free sheaves
on \m{X_2}, locally isomorphic to \
\m{(\ko_{X_2}\ot M)\oplus(\ko_X\ot N)}.

Let $\kh_{m,n}$ be the sheaf of algebraic groups on $X$, such that for every
open subset \m{U\subset X}, \m{\kh_{m,n}(U)} is the group of automorphisms
$\phi$ of the \m{\ko_X}-module \ \m{\ko_X(U)\ot(M\oplus N)} \ such that \
\m{\phi(\ko_X\ot N)\subset\ko_X\ot N}. The cohomology set \m{H^1(X,\kh_{m,n})}
is the set of isomorphism classes of surjective morphisms \ \m{\E\to E} \ of
vector bundles on $X$, where $\E$, $E$ are of rank \m{m+n}, $n$ respectively.

Let \ \m{\phi\in\kg_{m,n}(U)}. Then \m{\phi_{|U}} is an automorphism of \
\m{\ko_X(U)\ot(M\oplus N)} \ and \Nligne
\m{\phi(\ko_X\ot N)\subset\ko_X\ot N} (\ie \
\m{\phi_{|U}\in\kh_{m,n}(U)}), and
\[\xymatrix@R=2pt{\rho(U):\kg_{m,n}(U)\ar[r] & \kh_{m,n}(U)\\ \phi\fmaps[r] &
\phi_{|U}} \]
is a morphism of groups. In this way we define a surjective morphism of
sheaves of algebraic groups, \m{\rho:\kg_{m,n}\to\kh_{m,n}}. For every \
\m{\begin{pmatrix}
A & B\\ C & E\end{pmatrix}\in\kg_{m,n}(U)}, we have
\[\rho(U)\Big(\begin{pmatrix} A & B\\ C & E\end{pmatrix}\Big) \ = \
\begin{pmatrix}A_0 & 0\\ C & E\end{pmatrix} \ , \]
where \m{A_0} is the automorphism of \m{\ko_U\ot M} induced by $A$.
Let \ \m{\Gamma_{m,n}=\ker(\rho)}.

\sepprop

\begin{subsub}\label{lem1b}{\bf Lemma: } We have \
\m{\Gamma_{m,n}=L\ot\kl(M\oplus N,M)} \ (sheaf of additive groups).
\end{subsub}
\begin{proof}
$\ker(\rho(U))$ is the set of matrices \m{\dsp\begin{pmatrix}I_M+zC & zD \\
0 & I_N\end{pmatrix}}, with \ \m{C\in\ko_X(U)\ot\End(M)},
\m{D\in\ko_X(U)\ot \kl(N,M)}, and we have
\[\begin{pmatrix}I_M+zC & zD \\0 & I_N\end{pmatrix}
\begin{pmatrix}I_M+zC' & zD' \\0 & I_N\end{pmatrix} \ = \
\begin{pmatrix}I_M+z(C+C') & z(D+D') \\0 & I_N\end{pmatrix} \ . \]
The result follows immediately.
\end{proof}

\sepprop

Since $\Gamma_{m,n}$ is a sheaf of abelian groups, we can use the results of
\ref{coh_gr3}.

\sepprop

\begin{subsub}\label{prop1b}{\bf Proposition: } Let \ \m{g\in
H^1(X,\kh_{m,n})},
corresponding to the surjective morphism \ \m{\E\to E}. Then \
\m{\Gamma_{m,n}^g\simeq\E^*\ot E\ot L}.
\end{subsub}
\begin{proof}
let \m{(U_i)_{i\in I}} be an affine open cover of $X$, such that for every
\m{i\in I}, \m{L_{|U_i}} is trivial. Let \m{z_i} be a generator of
\m{L_{|U_i}}.
We have \ \m{z_j=\alpha_{ij}z_i}, where \ \m{\alpha_{ij}\in\ko_X(U_{ij})^*},
and \m{(\alpha_{ij})} is a cocycle that defines $L$. Let \m{(g_{ij})} be a
cocycle that represents $g$. We can write \m{g_{ij}} as a matrix \
\m{\dsp\begin{pmatrix}A_{ij} & 0\\ C_{ij} & E_{ij}\end{pmatrix}}, where \
\m{A_{ij}\in\Aut(\ko_X(U_{ij})\ot M)}, \m{E_{ij}\in\Aut(\ko_X(U_{ij})\ot N)} \
and \ \m{C_{ij}\in\ko_X(U_{ij})\ot\kl(M,N)} .
The cocycle \m{(g_{ij})} also represents $\E$, and \m{(A_{ij})} (resp.
\m{(\alpha_{ij}A_{ij})}) represents $E$ (resp. \m{E\ot L}). The vector bundle
\m{\E^*\ot E\ot L} is represented by \m{(\gamma_{ij})}, where
\[\xymatrix@R=2pt{\gamma_{ij}:\kl(M\oplus N,M)\ar[r] & \kl(M\oplus N,M) \\
\psi\fmaps[r] & \alpha_{ij}A_{ij}\circ\psi\circ g_{ij}^{-1} \ . }\]
Let \ \m{\psi=(P,Q)}, with \ \m{P\in \End(M)}, \m{Q\in \kl(N,M)}.
Then
\[\gamma_{ij}(\psi) \ = \ \Big(\alpha_{ij}\big(A_{ij}PA_{ij}^{-1}-
A_{ij}QE_{ij}^{-1}C_{ij}A_{ij}^{-1}\big),\alpha_{ij}A_{ij}QE_{ij}^{-1}\Big) \ .
\]
On the other hand, we have
\[g_{ij}\begin{pmatrix}I_m+z_jP & z_jQ\\ 0 & I_n\end{pmatrix}g_{ij}^{-1} \ = \
\begin{pmatrix}I_m+z_iP' & z_iQ'\\ 0 & I_n\end{pmatrix} \ , \]
with
\[P' \ = \ \alpha_{ij}\big(A_{ij}PA_{ij}^{-1}-
A_{ij}QE_{ij}^{-1}C_{ij}A_{ij}^{-1}\big) \ , \qquad Q' \ = \
\alpha_{ij}A_{ij}QE_{ij}^{-1} \ . \]
The result follows immediately.
\end{proof}

\sepprop

\begin{subsub}\label{obs_exi} Obstruction to the existence of quasi locally
free sheaves --  \rm Let \ \m{p:\E\to E} be a surjective morphism of vector
bundles on $X$, with \ \m{\rk{\E}=m+n}, \m{\rk(E)=m}.

According to proposition \ref{prop1b} and \ref{coh_gr}, {\em one can associate
to $p$ an element \Nligne \m{\nabla'_{\E,E}(p)\in\Ext^2_X(\E,E\ot L)} such that
there exists a quasi locally free sheaf $\ke$ on \m{X_2} such that \
\m{p\simeq\BS{\Psi}_\ke} \ if and only if \ \m{\nabla'_{\E,E}(p)=0}.} In
particular, this is the case if \ \m{\Ext^2_X(\E,E\ot L)=\nsp}.
\end{subsub}

\sepprop

\begin{subsub}\label{loc_f} Sheaves of groups associated to vector
bundles --
\rm Vector bundles are particular cases of quasi locally free sheaves. Hence we
can apply the preceding results to the extension of vector bundles on $X$ to
\m{X_2}. If $V$ is a vector bundle on $X$, the corresponding surjective
morphism is \ \m{I_V:V\to V}. Using the construction of $(\ref{equ5})$ by \v
Cech cohomology it is easy to see that this element associated to $V$ and that
of \ref{ext_hig} are the same, \ie
\m{\nabla'_{V,V}(I_V)=\nabla(V)=\xi_{X_2}(\nabla_0(V))}.

The sheaf of groups associated to $\E$ (resp. $E$) is \m{\kg_{m+n,m+n}} (resp.
\m{\kg_{m,m}}). Let
\[\eta_p:\Ext^2_X(\E,\E\ot L)\lra\Ext^2_X(\E,E\ot L) \ , \qquad
\zeta_p:\Ext^2_X(E,E\ot L)\lra\Ext^2_X(\E,E\ot L)\]
be the maps induced by $p$. As in the proof of proposition \ref{prop1b}, an
element \ \m{\gamma\in\Ext^2_X(\E,\E\ot L)} is represented by a cocycle
\m{(M_{ijk})}, where \m{M_{ijk}} is a matrix
\[M_{ijk} \ = \ \begin{pmatrix}\BS{A}_{ijk} & \BS{B}_{ijk} \\ \BS{C}_{ijk} &
\BS{E}_{ijk}\end{pmatrix} \
, \]
with \ \m{\BS{A}_{ijk}\in\ko_{U_{ijk}}\ot\kl(M,M)},
\m{\BS{B}_{ijk}\in\ko_{U_{ijk}}\ot\kl(N,M)},
\m{\BS{C}_{ijk}\in\ko_{U_{ijk}}\ot\kl(M,N)},\Nligne
\m{\BS{E}_{ijk}\in\ko_{U_{ijk}}\ot\kl(N,N)}. In the isomorphism of
\ref{prop1b}, applied to \ \m{I_\E:\E\to\E},
this corresponds to
\[\begin{pmatrix}I_m+z_i\BS{A}_{ijk} & z_i\BS{B}_{ijk} \\ z_i\BS{C}_{ijk} &
I_n+ z_i\BS{E}_{ijk}\end{pmatrix} \ \in \ \kg_{m+n,m+n}(U_{ijk}) \ . \]
Then \m{\eta_p(\gamma)} is represented by the cocycle
\m{(\BS{A}_{ijk},\BS{B}_{ijk})}, corresponding to
\[\begin{pmatrix}I_m+z_i\BS{A}_{ijk} & z_i\BS{B}_{ijk} \\ 0 & I_n\end{pmatrix}
\ \in \ \kg_{m,n}(U_{ijk}) \ . \]
Similarly, an element \ \m{\mu\in\Ext^2_X(E,E\ot L)} is represented by a
cocycle \m{(\BS{G}_{ijk})}, where \Nligne
\m{\BS{G}_{ijk}\in\ko_{U_{ijk}}\ot\kl(M,M)},
and \m{\zeta_p(\mu)} is represented by the cocycle \m{(\BS{G}_{ijk},0)},
corresponding to
\[\begin{pmatrix}I_m+z_i\BS{G}_{ijk} & 0 \\ 0 & I_n\end{pmatrix} \ \in \
\kg_{m,n}(U_{ijk}) \ . \]
\end{subsub}

\sepprop

\begin{subsub}\label{theo3b}{\bf Theorem:} We have \
\m{\nabla'_{\E,E}(p)=\eta_p(\nabla(\E))=\zeta_p(\nabla(E))} .
\end{subsub}\begin{proof} As in the proof of proposition \ref{prop1b}, let
\m{(g_{ij})} be a cocycle that represents the element of \m{H^1(X,\kh_{m,n})}
corresponding to $p$. We can write \m{g_{ij}} as a matrix \
\m{\dsp\begin{pmatrix}\BS{A}_{ij} & 0\\ \BS{C}_{ij} &
\BS{E}_{ij}\end{pmatrix}}, where \ \m{\BS{A}_{ij}\in\Aut(\ko_X(U_{ij})\ot M)},
\m{\BS{E}_{ij}\in\Aut(\ko_X(U_{ij})\ot N)} \
and \ \m{\BS{C}_{ij}\in\ko_X(U_{ij})\ot\kl(M,N)}.
The cocycle \m{(g_{ij})} also represents $\E$, and \m{(\BS{E}_{ij})} (resp.
\m{(\alpha_{ij}\BS{E}_{ij})}) represents $E$ (resp. \m{E\ot L}). Let \ \m{{\bf
g}_{ij}\in\kg_{m+n,m+n}(U_{ij})} \ be over \m{g_{ij}}. Suppose that \ \m{{\bf
g}_{ij}^{-1}={\bf g}_{ji}} \ for every distinct \m{i,j\in I}. Then \m{{\bf
g}_{ij}{\bf g}_{jk}{\bf g}_{ki}\in\Gamma_{m+n,m+n}(U_{ijk})}, and \m{({\bf
g}_{ij}{\bf g}_{jk}{\bf g}_{ki})} is a cocycle of \
\m{\Gamma_{m+n,m+n}^g=\Ext^2_X(\E,\E\ot L)} \ that represents
\m{\nabla(\E)} (cf. \ref{coh_gr}). Suppose that \ \m{{\bf
g}_{ij}(\ko_{X_2}(U_{ij}^{(2)})\ot N)\subset\ko_X(U_{ij})\ot N} \ for every
\m{i,j\in I}. Then \m{g_{ij}} induces \ \m{\ov{{\bf
g}}_{ij}\in\kg_{m,n}(U_{ij})}, which is also over \m{g_{ij}}, and such that \
\m{\ov{\bf g}_{ij}^{-1}=\ov{\bf g}_{ji}}. It follows that \m{(\ov{\bf
g}_{ij}\ov{\bf g}_{jk}\ov{\bf g}_{ki})} is a cocycle of \
\m{\Gamma_{m,n}^g=\Ext^2_X(\E,E\ot L)} \ that represents
\m{\nabla'_{\E,E}(p)}. The equalities of theorem \ref{theo3} follow from
\ref{loc_f} (using the fact that \m{\BS{B}_{ijk}=0} in \ref{loc_f}).

We now prove that such a family \m{({\bf g}_{ij})} exists.
Suppose that for every \m{i\in I}, there is a trivialization
\xmat{\delta_i:U_i^{(2)}\ar[r]^-\simeq & U_i\times{\bf Z}_2 \ . }
Then \ \m{\delta_{ij}^*=\delta_i^{*-1}\delta_j^*} \ is an automorphism of \
\m{\ko_{U_i\times Z_2}(U_i\times Z_2)=\ko_X(U_{ij})[t]/(t^2)}, such that for
every \ \m{\phi\in\ko_X(U_{ij})},  the term of degree zero of
\m{\delta_{ij}^*(\phi)} is the same as the term of degree zero of $\phi$. We
have the following description of \m{\delta_{ij}^*}:
\begin{enumerate}
\item[--] there is a derivation \m{D_{ij}} of \m{\ko_X(U_{ij})} such that \
\m{\delta_{ij}^*(\beta)=\beta+D_{ij}(\beta)t} \ for every
\m{\beta\in\ko_X(U_{ij})}.
\item[--] $\delta_{ij}^*(t)=\alpha_{ij}t$ .
\end{enumerate}
(cf. \ref{PMS}).  Since \
\m{U_{ij}\times{\bf Z}_2=\spec(\ko_{U_{ij}}[t]/(t^2))}, we can see \m{g_{ij}}
as an automorphism of \m{\ko_{U_{ij}\times{\bf Z}_2}\ot(M\oplus N)}. Let \
\m{\Phi_{ij}=\delta_i^*(g_{ij})}. It is an automorphism of
\m{\ko_{U_{ij}^{(2)}}\ot(M\oplus N)} over \m{g_{ij}} such that \
\m{\Phi_{ij}(\ko_{U_{ij}^{(2)}}\ot N)\subset\ko_{U_{ij}^{(2)}}\ot N}. But
we may not have \ \m{\Phi_{ij}^{-1}=\Phi_{ji}}. So instead of \m{g_{ij}} we take
\[h_{ij} \ = \ g_{ij}+\frac{1}{2}D_{ij}(g_{ij})t \ \ \in \
\Aut(\ko_{U_{ij}\times{\bf Z}_2}\ot\C^{m+n})\]
(it is a particular case of \cite{dr10}, 2.3), and \ \m{{\bf
g}_{ij}=\delta_i^*(h_{ij})}. We have \Nligne \m{{\bf
g}_{ij}(\ko_{X_2}(U_{ij}^{(2)})\ot\C^n)\subset\ko_X(U_{ij})\ot\C^n}. The
equality \ \m{{\bf g}_{ij}^{-1}={\bf g}_{ji}} \ is equivalent to \
\m{h_{ij}^{-1}=\delta_{ij}^*(h_{ji})}. We have \ \m{\dsp
D_{ji}=-\frac{1}{\alpha_{ij}}D_{ij}} \ and \
\m{D_{ij}(g_{ij}^{-1})=-g_{ij}^{-1}D_{ij}(g_{ij})g_{ij}^{-1}}. Hence
\begin{eqnarray*}
h_{ji} & = & g_{ij}^{-1}+\frac{1}{2}D_{ji}(g_{ij}^{-1})t \\
& = & g_{ij}^{-1}+\frac{1}{2\alpha_{ij}}g_{ij}^{-1}D_{ij}(g_{ij})g_{ij}^{-1}t \
,
\end{eqnarray*}
\begin{eqnarray*}
\delta_{ij}^*(h_{ji}) & = & g_{ij}^{-1}+\big(D_{ij}(g_{ij}^{-1})+
\frac{1}{2}g_{ij}^{-1}D_{ij}(g_{ij})g_{ij}^{-1}\big) \\
& = & g_{ij}^{-1}-\frac{1}{2}g_{ij}^{-1}D_{ij}(g_{ij})g_{ij}^{-1}t \ .
\end{eqnarray*}
On the other hand
\begin{eqnarray*}
h_{ij}^{-1} & = & \big(g_{ij}+\frac{1}{2}D_{ij}(g_{ij})t\big)^{-1}\\
& = & \Big[g_{ij}\big(1+g_{ij}^{-1}D_{ij}(g_{ij})\big)\Big]^{-1}\\
& = & g_{ij}^{-1}-\frac{1}{2}g_{ij}^{-1}D_{ij}(g_{ij})g_{ij}^{-1}t \ .
\end{eqnarray*}
\end{proof}

\sepprop

\begin{subsub}\label{coro4}{\bf Corollary:} We have \
\m{\nabla'_{\E,E}(p)=\nabla_{\E,E}(p)}.
\end{subsub}
\begin{proof}
This follows from theorems \ref{theo3}, \ref{theo3b}.
\end{proof}

\end{sub}

\sepsub

\Ssect{Exact sequences of length 4}{L4}

Let $\ke$ be a coherent sheaf on \m{X_2}. Recall that \ \m{\BS{\Psi}_\ke} \
denotes the morphism \ \m{\ke_{|X}\to\ke_1\ot L^*} \ induced by \
\m{\ke\ot\ki_{X,X_2}\to\ke} (induced by \ \m{\ki_{X,X_2}\subset\ko_{X_2}}). We
have a canonical exact sequence of sheaves on $X$
\xmat{\Sigma(\ke): \ 0\ar[r] & \ke_1\ar[r]^-\iota & \ke^{(1)}\ar[r] &
\ke_{|X}\ar[r]^-{\BS{\Psi}_\ke} & \ke_1\ot L^*\ar[r] & 0\ , }
where $\iota$ is the inclusion. If $\ke$ is quasi locally free it is an exact
sequence of vector bundles. If $\ke$ locally isomorphic to \
\m{(\ko_{X_2}\ot\C^m)\oplus(\ko_X\ot\C^n)}, \m{\ke_1} is of rank $m$, and
\m{\ke^{(1)}}, \m{\ke_{|X}} are of rank \m{m+n}.

We will be interested in exact sequences of vector bundles on open subsets $U$
of $X$
\[0\lra E\ot L_{|U}\lra\F\lra\G\lra E\lra 0 \ , \]
where $E$ is of rank $m$, and $\E$, $\F$ of rank \m{m+n}. Given two such exact
sequences
\[\Sigma: \ 0\lra E\ot L_{|U}\lra\F\lra\G\lra E\lra 0 \ , \]
\[\Sigma': \ 0\lra E'\ot L_{|U}\lra\F'\lra\G'\lra E'\lra 0 \ , \]
an {\em isomorphism} \m{\Sigma\to\Sigma'} is a triple
\m{(\epsilon,\phi,\gamma)} of isomorphisms \ \m{\epsilon:E\to E'},
\m{\phi:\F\to\F'}, \m{\gamma:\G\to\G'} such that the following diagram is
commutative
\xmat{\Sigma: \ 0\ar[r] & E\ot L_{|U}\ar[r]\ar[d]^{\epsilon\ot I_L} &
\F\ar[r]\ar[d]^\phi & \G\ar[r]\ar[d]^\gamma & E\ar[r]\ar[d]^\epsilon & 0 \\
\Sigma': \ 0\ar[r] & E'\ot L_{|U}\ar[r] & \F'\ar[r] & \G'\ar[r] & E'\ar[r] & 0}
If \m{V\subset U} is an open subset, the restriction \m{\Sigma_{|V}} of
$\Sigma$ to $V$ is defined in an obvious way.

Let \m{\ke}, \m{\ke'} be quasi locally free sheaves on \m{U^{(2)}}, locally
isomorphic to \ \m{(\ko_{U^{(2)}}\ot\C^m)\oplus(\ko_U\ot\C^n)}. an isomorphism
\ \m{f:\ke\to\ke'} \ induces an obvious isomorphism \
\m{\Sigma(f):\Sigma(\ke)\to\Sigma(\ke')}.

\sepprop

\begin{subsub}\label{as_K} The associated sheaf of algebraic groups -- \rm Let
$M$, $N$ be complex vector spaces, with \ \m{\dim(M)=m}, \m{\dim(N)=n}.
The  {\em canonical exact sequence on $U$} is
\[C_U: \ 0\lra L_{|U}\ot M\hfl{\iota}{}(L_{|U}\ot M)\oplus(\ko_U\ot N)
\hfl{\zeta}{}\qquad\qquad \]
\[\qquad\qquad(\ko_U\ot M)\oplus(\ko_U\ot N)\hfl{\nu}{}\ko_U\ot M\lra 0 \ , \]
with
\[\iota \ = \ \begin{pmatrix}I\\ 0\end{pmatrix} \ , \quad \xi \ = \
\begin{pmatrix}0 & 0\\ 0 & I\end{pmatrix} \ , \quad \nu \ = \ \begin{pmatrix}I
& 0\end{pmatrix} \ . \]
i.e \ \m{\Sigma_U=\Sigma((\ko_{U^{(2)}}\ot M)\oplus(\ko_U\ot N))}.

If $\Sigma$ is an exact sequence on $U$, and \m{x\in U} is a closed point,
there exists an open neighborhood $V$ of $x$ in $U$ such that \
\m{\Sigma_{|V}\simeq C_V}.

Let $\kk_{m,n}$ the sheaf of algebraic groups on $X$, such that for every open
subset \m{U\subset X}, \m{\kk_{m,n}(U)} is the group of automorphisms of \
\m{C_U}, \ie it is the group of triples
\[\begin{pmatrix}
A_0\ot I_L \ , \quad \begin{pmatrix}A_0\ot I_L & B\\ 0 & E\end{pmatrix} \ ,
\quad\begin{pmatrix}A_0 & 0 \\ C & E\end{pmatrix}\end{pmatrix} \ , \]
where \m{A_0} (resp . $E$) is an automorphism of \ \m{\ko_U\ot M} (resp. \
\m{\ko_U\ot N}),\Nligne \m{B\in\Hom(\ko_U\ot N, L\ot M)}, \m{C\in\Hom(\ko_U\ot
M,\ko_U\ot N)}. It can be defined as well as the group of quadruples
\m{(A_0,B,C,E)}, with the group law
\[(A_0,B,C,E)*(A'_0,B',C',E') \ = \ (A_0A'_0,(A\ot I_L)B'+BE',
CA'_0+EC'_0,EE') \ . \]
Let
\[\xymatrix@R=5pt{ \Theta(U):\kg_{m,n}(U)\ar[r] & \kk_{m,n}(U)\\
\phi\fmaps[r] & \Sigma(\phi)} \]
We define in this way a surjective morphism of sheaves of groups \
\m{\Theta:\kg\to\kk}. We have
\[\Theta(U)\Big(\begin{pmatrix}A & B\\C & E\end{pmatrix}\Big) \ = \
(A_0,B,C,E)
\ , \]
where \m{A_0} is the automorphism of \ \m{\ko_U\ot M} \ induced by $A$. Let \
\m{\Xi_{m,n}=\ker(\Theta)}, which isomorphic to \m{L\ot\End(M)}.
\end{subsub}

\sepprop

\begin{subsub}\label{prop1c}{\bf Proposition: } Let \ \m{h\in
H^1(X,\kk_{m,n})}, corresponding to
\[\Sigma: \ 0\lra E\ot L_{|U}\lra\F\lra\G\lra E\lra 0 \ . \]
Then \ \m{\Xi_{m,n}^h\simeq E^*\ot E\ot L}.
\end{subsub}
\begin{proof}
Similar to that of proposition \ref{prop1b}.
\end{proof}

\sepprop

According to proposition \ref{prop1c} and \ref{coh_gr}, {\em one can associate
to $\Sigma$ an element \Nligne \m{\nabla''(\Sigma)\in\Ext^2_X(E,E\ot L)} such
that there exists a quasi locally free sheaf $\ke$ on \m{X_2} such that \
\m{\Sigma\simeq\Sigma(\ke)} \ if and only if \ \m{\nabla''(\Sigma)=0}.} In
particular, this is the case if \ \m{\Ext^2_X(E,E\ot L)=\nsp}.

We have a canonical morphism \ \m{\Upsilon(U):\kk_{m,n}(U)\to\kh_{m,n}(U)},
such that
\[\Upsilon(U)\Big(\begin{pmatrix}
A_0\ot I_L \ ,\begin{pmatrix}A_0\ot I_L & B\\ 0 & E\end{pmatrix} \ ,
\begin{pmatrix}A_0 & 0 \\ C & E\end{pmatrix}\end{pmatrix}\Big) \ = \
\begin{pmatrix}A_0 & 0 \\ C & E\end{pmatrix} \ . \]

We define in this way a morphism of sheaves of groups \
\m{\Upsilon:\kk_{m,n}\to\kh_{m,n}}, and we have a commutative diagram
\xmat{\Xi_{m,n}\flinc[r]\ar[d]^\eta & \kg_{m,n}\fleq[d]\ar[rr]^-\Theta & &
\kk_{m,n}\ar[d]^\Upsilon\\
\Gamma_{m,n}\flinc[r] & \kg_{m,n}\ar[rr]^-\rho & & \kh_{m,n}}
Let \ \m{h\in H^1(X,\kk_{m,n})} \ and \ \m{g=H^1(\Upsilon)(h)}. Suppose that
$g$ corresponds to the exact sequence
\[\Sigma:0\lra E\ot L_{|U}\lra\F\lra\E\hfl{p}{} E\lra 0 \ . \]
Then \m{H^1(\Upsilon)(h)} corresponds to \ \m{p:\E\to E}.

We have a canonical induced morphism
\[\eta_g:\Xi_{m,n}^h=E^*\ot E\ot L\lra\Gamma_{m,n}^g=\E^*\ot E\ot L \ . \]

\sepprop

let \ \m{N=\ker(p)}. For every \ \m{\epsilon\in\Ext^1_X(N,E\ot L)},
corresponding to an exact sequence \Nligne \m{0\to E\ot
L\hfl{\iota_\epsilon}{}\F\to N\to 0}, let \m{\Sigma_\epsilon} be the exact
sequence
\[0\lra E\ot L\hfl{\iota_\sigma}{}\F\to\E\hfl{p}{} E\to 0 \ . \]

We have an exact sequence
\[\Ext^1_X(\E,E\ot L)\lra\Ext^1_X(N,E\ot L)\hfl{\partial_p}{}
\Ext^2_X(E,E\ot L)\hfl{\mu_p}{}\Ext^2_X(\E,E\ot L)\]
Using the methods of \ref{const_QLF}, one proves easily

\sepprop

\begin{subsub}\label{prop8}{\bf Proposition: } {\bf 1 -- } \m{\eta_g} is the
canonical morphism induced by $p$.

{\bf 2 -- } We have \ \m{\nabla_{\E,E}(p)=\mu_p(\nabla''(\Sigma))} .

{\bf 3 -- } Let \ \m{\epsilon_0,\epsilon_1\in\ker(\mu_p)}. Then we have
\[\nabla''(\Sigma_{\epsilon_0})-\nabla''(\Sigma_{\epsilon_1}) \ = \
\partial_p(\epsilon_0-\epsilon_1) \ . \]
\end{subsub}

\end{sub}

\sepsec

\section{Examples -- the case of \texorpdfstring{\m{\P_2}}{5}}\label{ex-p2}

There are two primitive double schemes with reduced subscheme \m{\P_2} and
associated line bundle \ \m{L=\omega_{\P_2}=\ko_{\P_2}(-3)} : the trivial one
\m{\BS{X}_2^0}, and the non-projective one \m{\BS{X}_2}.

We have \ \m{\Pic(\BS{X}_2^0)\simeq\Z} : for every \m{k\in\Z}, there is exactly
one extension of \m{\ko_{\P_2}(k)} to a line bundle \m{\ko_{\BS{X}_2^0}(k)} on
\m{\BS{X}_2^0}. If \ \m{\pi:\BS{X}_2^0\to\P_2} \ is the canonical projection,
we have \ \m{\ko_{\BS{X}_2^0}(k)=\pi^*(\ko_{\P_2}(k))}.

The only line bundle on \m{\BS{X}_2} is \m{\ko_{\BS{X_2}}}.

Let \m{X_2} be one of the primitive multiple schemes
\m{\BS{X}_2}, \m{\BS{X}_2^0}.

Let $\E$, $E$ be vector bundles on \m{\P_2}. Recall that we have a canonical
exact sequence
\[0\lra\Ext^1_{\P_2}(\E,E(-3))\lra\Ext^1_{X_2}(\E,E(-3))
\hfl{\gamma_{\E,E}}{}\Hom(\E,E)\qquad\qquad\qquad\]
\[\qquad\qquad\qquad\qquad\qquad \hfl{\nabla_{\E,E}}{}\Ext^2_{\P_2}(\E,E(-3))
\]
Let $\ke$ be a quasi locally free sheaf on \m{X_2}, \m{\E=\ke_{|\P_2}},
\m{E=\ke_1(3)}, and \ \m{\sigma\in\Ext^1_{X_2}(\E,E(-3))} \ corresponding to
the canonical extension \ \m{0\to E(-3)\to\ke\to\E\to 0}.
Then \ \m{\gamma_{\E,E}(\sigma)} is the canonical surjective morphism \
\m{\BS{\Psi}_\ke:\ke_{|\P_2}\lra\ke_1(3)} .

\sepsub

\Ssect{Rank 3}{PR3}

We will study quasi locally free sheaves which are locally isomorphic to \
\m{\ko_{X_2}\oplus\ko_{\P_2}}. If $\ke$ is such a sheaf, We have \
\m{\rk(\ke_1(3))=1} \ and \
\m{\rk(\ke_{|X})=2}, hence \ \m{\rk(\ker(\BS{\Psi}_\ke))=1}.

\sepprop

We have \ \m{\ke_1\ot L^*=\ke_1(3)\simeq\ko_{\P_2}(m)},
\m{\ker(\BS{\Psi}_\ke)\simeq\ko_{\P_2}(n)} \
for some \ \m{m,n\in\Z}, and an exact sequence
\[0\lra\ko_{\P_2}(n)\lra\ke_{|\P_2}\lra\ko_{\P_2}(m)\lra 0 \ , \]
and since \ \m{\Ext^1_{\P_2}(\ko_{\P_2}(m),\ko_{\P_2}(n))=\nsp}, \m{\ke_{|\P_2}=
\ko_{\P_2}(m)\oplus\ko_{\P_2}(n)} \ and \m{\BS{\Psi}_\ke} is the projection \
\m{p:\ko_{\P_2}(m)\oplus\ko_{\P_2}(n)\to\ko_{\P_2}(m)}.

So we will suppose that \ \m{\E=\ko_{\P_2}(m)\oplus\ko_{\P_2}(n)} \ and
\m{E=\ko_{\P_2}(m)}, and \ \m{p:\E\to E} \ is the projection. The canonical map
\xmat{\mu_p :\Ext^2_X(E,E\ot L)\fleq[d]\ar[rr] & &
\Ext^2_X(\E,E\ot L)\fleq[d] \\
H^2(\P_2,\ko_{\P_2}(-3)) & & H^2(\P_2,\ko_{\P_2}(-3))\oplus
H^2(\P_2,\ko_{\P_2}(m-n-3))}
is the inclusion.

It follows that \ \m{\nabla_{\E,E}(p)=0} \ if and only if \ \m{\nabla(E)=0}. On
the other hand, since \ \m{\Ext^1_{\P_2}(\E,E\ot L)=\nsp}, a quasi locally free
sheaf such that \ \m{\BS{\Psi}_\ke=p} \ is unique if it exists. Hence we have

\sepprop

\begin{subsub}\label{prop10}{\bf Proposition: } The quasi locally free sheaves
on \m{X_2} locally isomorphic to \ \m{\ko_{X_2}\oplus\ko_X} \ are
\begin{enumerate}
\item[--] the sheaves \ \m{\ko_{X_2}(m)\oplus\ko_{\P_2}(n)}, \m{m,n\in\Z}, if \
\m{X_2=\BS{X}_2^0},
\item[--] the sheaves  \ \m{\ko_{X_2}\oplus\ko_{\P_2}(n)}, \m{n\in\Z}, if \
\m{X_2=\BS{X}_2}.
\end{enumerate}
\end{subsub}

\end{sub}

\sepsec

\section{Examples -- the case of
\texorpdfstring{\m{\P_1\times\P_1}}{5}}\label{ex-p1p1}

Let \m{V_1}, \m{V_2} two complex vector spaces of dimension 2, and \
\m{\BS{\kx}=\P(V_1)\times\P_1(V_2)}. Let \ \m{p_1:\BS{\kx}\to\P(V_1)},
\m{p_2:\BS{\kx}\to\P(V_2)} \ be the projections. For \m{m_1,m_2\in\Z}, let \
\[\ko(m_1,m_2) \ = \ p_1^*(\ko_{\P(V_1)}(m_1))\ot p_2^*(\ko_{\P(V_2)}(m_2)) \ .
\]
We have \ \m{\omega_\BS{\kx}\simeq\ko(-2,-2)} \ and \
\m{T_\BS{\kx}\simeq\ko(2,0)\oplus\ko(0,2)} .

We have an isomorphism \ \m{\Z\times\Z\simeq\Pic(\BS{\kx})}, associating
\m{\ko(m_1,m_2)} to \m{(m_1,m_2)}.

\sepprop

\begin{subsub}\label{lem1}{\bf Lemma: } {\bf 1 -- }We have canonical
isomorphisms
\begin{eqnarray*}
H^0(\BS{\kx},\ko(m_1,m_2)) & \simeq & S^{m_1}V_1^*\ot S^{m_2}V_2^* \ \
\text{if} \ \ m_1,m_2 \geq 0 , \\
& = & 0 \ \ \text{if} \ \ m_1<0 \ \text{or} \ \ m_2<0 \ .
\end{eqnarray*}
\begin{eqnarray*}
H^1(\BS{\kx},\ko(m_1,m_2)) & \simeq & S^{m_1}V_1^*\ot S^{-2-m_2}V_2 \ \
\text{if} \ \ m_1\geq 0, m_2 \leq -2 , \\
& \simeq & S^{-m_1-2}V_1\ot S^{m_2}V_2^* \ \ \text{if} \ \ m_1\leq -2, m_2\geq
0 \ , \\
& = & 0 \ \ \text{in all the other cases} \ .
\end{eqnarray*}
\begin{eqnarray*}
H^2(\BS{\kx},\ko(m_1,m_2)) & \simeq & S^{-m_1-2}V_1\ot S^{-m_2-2}V_2 \ \
\text{if} \ \ m_1,m_2 \leq -2 , \\
& = & 0 \ \ \text{if} \ \ m_1>-2 \ \text{or} \ \ m_2>-2 \ .
\end{eqnarray*}
\end{subsub}

\sepprop

Let \m{\BS{\kx}_2} a primitive double scheme such that \
\m{(\BS{\kx}_2)_{red}=\BS{\kx}}, and with associated line bundle \
\m{L=\omega_{\BS{\kx}}}.

\sepsub

\Ssect{Line bundles on double structures on \ \m{\P_1\times\P_1}}{p1p1}

According to \ref{param_DS}, the primitive double schemes with associated
smooth scheme \m{\BS{\kx}} and associated line bundle \m{\omega_{\BS{\kx}}}
are parameterized by \ \m{\P(H^1(\BS{\kx},\omega_{\BS{\kx}}))\cup\nsp}, and
\m{H^1(\BS{\kx}, \Omega_\BS{\kx})} is generated by \
\m{e_1=\nabla_0(\ko(1,0))} \ and \ \m{e_2=\nabla_0(\ko(0,1))}. For the
multiplication map
\[H^1(\BS{\kx},\Omega_\BS{\kx})\times H^1(\BS{\kx},\Omega_\BS{\kx})\lra
H^2(\BS{\kx},\omega_\BS{\kx})\simeq\C \ , \]
we have \ \m{e_1^2=e_2^2=0}, \m{e_1.e_2=1}.

We have \ \m{H^1(\BS{\kx},\omega_\BS{\kx})=\nsp}. Hence, from \ref{ext_ss}, if
a line bundle on $\BS{\kx}$ can be extended to \m{\BS{\kx}_2}, this extension
is unique.

Suppose that \
\m{\sigma(\BS{\kx}_2)=\C(\alpha_1e_1+\alpha_2e_2)\in
\P(H^1(\BS{\kx},\Omega_\BS{\kx}))\cup\nsp}. Then from \ref{ext_hig},
\m{\ko(m_1,m_2)} can be extended to a line bundle on \m{\BS{\kx}_2} if and
only if \ \m{\alpha_1m_2+\alpha_2m_1=0}.

\sepprop

\begin{subsub}\label{line_0}{\bf Proposition: }{\bf 1 -- } If
\m{\alpha_1,\alpha_2\not=0} and \m{\dsp\frac{\alpha_1}{\alpha_2}} is
irrational, the only line bundle on \m{\BS{\kx}_2} is \m{\ko_{\BS{\kx}_2}},
and \m{\BS{\kx}_2} is not quasi-projective.

{\bf 2 -- } If \m{\alpha_1,\alpha_2\not=0} and \
\m{\dsp\frac{\alpha_1}{\alpha_2}=\frac{n_1}{n_2}}, where \m{n_1,n_2} are
relatively prime integers, we have \ \m{\Pic(\BS{\kx}_2)\simeq\Z}, generated by
a line bundle whose restriction to $\BS{\kx}$ is \m{\ko(n_1,-n_2)}. Moreover,
\m{\BS{\kx}_2} is quasi-projective if and only if \ \m{\alpha_1\alpha_2<0}.

{\bf 3 -- } If \m{\alpha_1=0} we have \ \m{\Pic(\BS{\kx}_2)\simeq\Z}. It is
generated by a line bundle whose restriction to $\BS{\kx}$ is \m{\ko(0,1)}. If
\m{\alpha_2=0} we have \ \m{\Pic(\BS{\kx}_2)\simeq\Z}. It is generated by
a line bundle whose restriction to $\BS{\kx}$ is \m{\ko(1,0)}. In the two
cases, \m{\BS{\kx}_2} is not quasi-projective.
\end{subsub}
\begin{proof}
This follows easily from \ref{ext_hig}.
\end{proof}

\end{sub}

\sepsub

\Ssect{Rank 3}{PR3_b}

We will study quasi locally free sheaves which are locally isomorphic to \
\m{\ko_{\BS{\kx}_2}\oplus\ko_\BS{\kx}}. If $\ke$ is such a sheaf, we have \
\m{\rk(\ke_1(3))=1} \ and \ \m{\rk(\ke_{|\BS{\kx}})=2}, hence \
\m{\rk(\ker(\BS{\Psi}_\ke))=1}.
We have \
\m{\ke_1\ot L^*=\ke_1(2,2)=\ko(a_1,b_1)} , \m{\ker(\BS{\Psi}_\ke)=\ko(a_0,b_0)}
\ for some \m{a_0,b_0,a_1,b_1\in\Z}, and an exact sequence
\[0\lra\ker(\BS{\Psi_\ke})=\ko(a_0,b_0)\lra\ke_{|\BS{\kx}}
\hfl{\BS{\Psi_\ke}}{}\ke_1(2,2)=\ko(a_1,b_1)\lra 0 \ . \]
If the extension is trivial, then \m{\BS{\Psi_\ke}} is the projection \
\m{p:\ko(a_0,b_0)\oplus\ko(a_1,b_1)\to\ko(a_1,b_1)} . Conversely, consider the
projection $p$. Then as in \ref{PR3}, there exists a quasi locally free sheaf
$\ke$ on \m{\BS{\kx}_2} such that \ \m{\BS{\Psi}_\ke\simeq p} \ if and only if
\m{\ko(a_1,b_1)} can be extended to a line bundle on \m{\BS{\kx}_2}. These line
bundles are described in \ref{line_0}.

We now consider non trivial extensions
\begin{equation}\label{equ10}0\lra\ko(a_0,b_0)\lra\E\hfl{\Phi}{}
E=\ko(a_1,b_1)\lra 0 \ . \end{equation}
By lemma \ref{lem1}, they exist only if
\begin{equation}\label{equ12}
a_0-a_1 \ \geq 0 \ , \ b_0-b_1\leq -2 \qquad \text{or} \qquad
a_0-a_1 \ \leq \ -2 \ , \ b_0-b_1\geq 0 \ .\end{equation}
We suppose now that these conditions are satisfied.

If \ \m{a_0=a_1} \ or \ \m{b_0=b_1}, the extension $(\ref{equ10})$ is the
pullback
of an extension on one of the \m{\P_1}, twisted by a line bundle on
$\BS{\kx}$. If \ \m{a_0\not=a_1} \ and \ \m{b_0\not=b_1}, $\E$ is simple, by
lemma \ref{lem5}.

\sepprop

\begin{subsub}\label{prop20}{\bf Proposition: } There exist quasi locally free
sheaves $\ke$ on \m{\BS{\kx}_2} such that \ \m{\BS{\Psi}_\ke\simeq\Phi}.
\end{subsub}
\begin{proof}
We will prove that \ \m{\Ext^2_\BS{\kx}(\E,E\ot L)\simeq\Hom(E,\E)^*=
\nsp}. The result follows then from corollary \ref{coro2}.

Suppose first that \ \m{a_0\not=a_1} \ and \ \m{b_0\not=b_1}. Since $\E$ is
simple, by proposition \ref{prop4b} we have \ \m{\Ext^2_\BS{\kx}(\E,E\ot
L)=\nsp}.

Suppose that \ \m{a_0=a_1} (the case \ \m{b_0=b_1} \ is similar). Then we have
 \ \m{\E=p_2^*(\E')\ot\ko_(a_1,b_1)}, where \m{\E'} is a vector bundle on
\m{\P_1},
and there is a non trivial extension
\[0\lra\ko_{\P_1}(b_0-b_1)\lra\E'\lra\ko_{\P_1}\lra 0 \ . \]
Since every vector bundle on \m{\P_1} is a direct sum of line bundles, we have
 \ \m{\E'=\ko_{\P_1}(p)\oplus\ko_{\P_1}(q)}, with (since the extension is not
trivial), \m{b_0-b_1<p\leq q<0} \ and \ \m{p+q=b_0-b_1}. Hence \
\m{\E\simeq\ko(a_0,p+b_1)\oplus\ko(a_0,q+b_1)}, and \ \m{\Hom(E,\E)\simeq
H^0(\P_1,\ko_{\P_1}(p))\oplus H^0(\P_1,\ko_{\P_1}(q))=\nsp}, since \ \m{p,q<0}.
\end{proof}

\sepprop

\begin{subsub}\label{prop23}{\bf Proposition: } Let $\ke$ be a quasi locally
free sheaf on \m{\BS{\kx}_2} such that \ \m{\BS{\Psi}_\ke\simeq\Phi}. Then
$\ke$ is simple, \m{h^1(\EEnd(\ke))=2h^1(\ko(a_0-a_1,b_0-b_1))-2} \ and \
\m{H^2(\BS{\kx}_2,\EEnd(\ke))=\nsp}.
\end{subsub}
\begin{proof} From $(\ref{equ10})$ and $(\ref{equ12})$ we have \
\m{\Hom(\E,E\ot L)=\nsp}, and by proposition \ref{prop12}, $\ke$ is simple.

Now we prove that \ \m{H^2(\BS{\kx}_2,\EEnd(\ke))=\nsp}.
Let \ \m{\F=\ke^{(1)}/\ke_1} . From the exact sequence \
\m{0\to\ke_1\to\ke\to\ke_{|\BS{\kx}}\to 0} \ and theorem
\ref{theo5} follows the exact sequence
\xmat{0\ar[r] & \HHom(\ke,\ke_1)\ar[r]\fleq[d] & \EEnd(\ke)\ar[r]^-\mu &
\HHom(\ke,\ke_{|\BS{\kx}})\ar[r]\fleq[d] & \qquad\qquad \\
 & \HHom(\ke_{|\BS{\kx}},\ke_1) & & \EEnd(\ke_{|\BS{\kx}})}
\[\qquad\qquad\qquad
\xymatrix{\EExt^1_{\BS{\kx}_2}(\ke,\ke_1)\ar[r]^-\alpha\fleq[d] &
\EExt^1_{\BS{\kx}_2}(\ke,\ke)\fleq[d]\\
\HHom(\F\ot L,\ke_1) & \HHom(\F\ot L,\F)}\]
We have \m{\alpha=0}, since it is induced by the morphism
\[\beta:(\ke_1)^{(1)}/(\ke_1)_1\lra \ke^{(1)}/\ke_1\]
which comes from the inclusion \ \m{\ke_1\subset\ke}, and \m{\beta=0}. Let \
\m{K=\imm(\mu)}. We have two exact sequences
\begin{equation}\label{equ7}
0\lra\HHom(\ke_{|\BS{\kx}},\ke_1)\lra\EEnd(\ke)\lra K\lra 0 \ ,
\end{equation}
\begin{equation}\label{equ8}0\lra
K\lra\EEnd(\ke_{|\BS{\kx}})\hfl{\delta}{}\HHom(\F\ot
L,\ke_1)\lra 0 \ . \end{equation}

Using $(\ref{equ10})$, it is easy to see that \
\m{H^2(\BS{\kx}_2,\HHom(\ke_{|\BS{\kx}},\ke_1))=\nsp} \ and \
\m{H^2(\BS{\kx}_2,\EEnd(\ke_{|\BS{\kx}}))=\nsp}. Hence \
\m{H^2(\BS{\kx}_2,\EEnd(\ke))\simeq\coker(H^1(\delta))}, and it suffices to
prove that \m{H^1(\delta)} is surjective. Let
\[\lambda:\ke_{|\BS{\kx}}\lra\ke_1\ot L^*\ , \qquad
\tau:(\ke_{|\BS{\kx}})^*\to\F^* \ . \]
be the morphisms deduced from the canonical exact sequence \
\m{0\to\kn\to\ke_{|\BS{\kx}}\to\ke_1\ot L^*\to 0}. By proposition
\ref{prop22}, we have \ \m{H^1(\delta)=
H^1(\lambda\ot I_{\F^*})\circ H^1(I_{\ke_{|\BS{\kx}}}\ot\tau)}, hence it
suffices to prove that \m{H^1(\lambda\ot I_{\F^*})} and
\m{H^1(I_{\ke_{|\BS{\kx}}}\ot\tau)} are surjective. This follows from the fact
that \ \m{H^2(\BS{\kx}_2,\F\ot\F^*)=\nsp} \ and \
\m{H^2(\BS{\kx}_2,\ke_{|\BS{\kx}}\ot\ke_1\ot L^*)=\nsp} (which can be deduced
easily from $(\ref{equ10})$).

The fact that \ \m{h^1(\EEnd(\ke))=2h^1(\ko(a_0-a_1,b_0-b_1))-2} \ follows
easily from $(\ref{equ7})$, $(\ref{equ8})$, and the exact sequence \
\m{0\to\ko(a_1,b_1)\to\ke_{|\BS{\kx}}\to\ko(a_0,b_0)\to 0}.
\end{proof}

\sepprop

\begin{subsub}{\bf Remark: }\rm Proposition \ref{prop23} can also be proved
by using the canonical exact sequence
\begin{equation}\label{equ11}0\lra\HHom(\ke_{|\BS{\kx}},\ke^{(1)})\lra
\EEnd(\ke)\lra\EEnd(\ke_1)\lra 0 \ . \end{equation}
\end{subsub}

\sepprop

\begin{subsub} Deformation of sheaves -- \rm By theorem \ref{theo5}, we have
\[\EExt^1_{\BS{\kx}_2}(\ke,\ke)\simeq\Hom((\ke^{(1)}/\ke_1)\ot
L,\ke^{(1)}/\ke_1)=\Hom(\ko(a_0,b_0)\ot L,\ko(a_0,b_0))\simeq\ko(2,2) \ . \]
By proposition \ref{prop23} and the Ext spectral sequence, we have an exact
sequence
\[0\lra H^1(\BS{\kx_2},\EEnd(\ke))\lra \Ext^1_{\BS{\kx}_2}(\ke,\ke)\lra
H^0(\BS{\kx},\ko(2,2))\lra 0 \ . \]
Let \ \m{\Sigma(a_0,b_0,a_1,b_1)} \ be the set of isomorphism classes of quasi
locally free sheaves $\ke$ on \m{\BS{\kx}_2} such that \
\m{\ke_1\simeq\ko(a_0-2,b_0-2)} \ and \ \m{\ke^{(1)}/\ke_1\simeq\ko(a_1,b_1)}.
This set is open (by theorem \ref{theo6}).

It follows from corollary \ref{coco5} that for every flat family $\kf$ of
sheaves in \m{\Sigma(a_0,b_0,a_1,b_1)} parameterized by a smooth scheme
$S$, and any closed point \m{s\in S}, the image of the the Koda\"ira-Spencer
map \ \m{T_{\pmb{\M},m}\to\Ext^1_{\BS{\kx}_2}(\kf_s,\kf_s)} \ is contained if
\m{H^1(\BS{\kx_2},\EEnd(\kf_s))}.
\end{subsub}

\sepprop

\begin{subsub}\label{fin_modPP} Fine moduli spaces -- \rm (cf. \ref{fam_qlf}).
We suppose that \ \m{a_0\not=a_1} \ and \ \m{b_0\not=b_1}. Let \Nligne
\m{H=\Ext^1_\BS{\kx}(\ko(a_1,b_1),\ko(a_0,b_0))}. Then if we have two exact
sequences on $\BS{\kx}$
\[0\to\ko(a_0,b_0)\to\E\to\ko(a_1,b_1)\to 0 \ , \qquad
0\to\ko(a_0,b_0)\to\E'\to\ko(a_1,b_1)\to 0 \ , \]
corresponding to \ \m{\sigma,\sigma'\in H} \ respectively, $\E$ and \m{\E'}
are isomorphic if and only if there exists \ \m{\lambda\in\C^*} \ such that
\ \m{\sigma'=\lambda\sigma}.

Let \ \m{\pmb{\P}=\P(H)}, which parameterizes the non trivial extensions
$(\ref{equ10})$, and \ \m{p_\BS{\kx}:\BS{\kx}\times\pmb{\P}\to\BS{\kx}}, \
\m{p_{\pmb{\P}}:\BS{\kx}\times\pmb{\P}\to\pmb{\P}} \ the projections. There
exists a universal extension
\[0\lra p_\BS{\kx}^*(\ko(a_0,b_0))\ot p_{\pmb{\P}}^*(\ko_{\pmb{\P}}(1))\lra
\pmb{\E}\lra p_\BS{\kx}^*(\ko(a_1,b_1))\lra 0 \]
on \ \m{\BS{\kx}\times\pmb{\P}}. We will see also $\pmb{\E}$ as a sheaf on \
\m{\BS{\kx}_2\times\pmb{\P}}.

We will now use \ref{fam_qlf} to build a fine moduli space for families of
sheaves in \m{\Sigma(a_0,b_0,a_1,b_1)} parameterized by integral schemes (cf.
\ref{fin-mod}).

Let \ \m{q_{\BS{\kx}_2}:\BS{\kx}_2\times\pmb{\P}\to\BS{\kx}_2}, \
\m{q_{\pmb{\P}}:\BS{\kx}_2\times\pmb{\P}\to\pmb{\P}} \ be the projections.

Let \ \m{0\to\ko(a_0,b_0)\to\E\to\ko(a_1,b_1)\to 0} \ be a non trivial
extension. Then we have \Nligne \m{\dim(\Hom(\E,\ko(a_1,b_1)))=1},
\m{\dim(\Ext^1_\BS{\kx}(\E,\ko(a_1-2,b_1-2)))=h^1(\ko(a_0-a_1,b_0-b_1))-1} \
and \ \m{\Ext^2_\BS{\kx}(\E,\ko(a_1-2,b_1-2))=\nsp}. It follows from theorem
\ref{theo1b} that \Nligne
\m{\Ext^1_{q_{\pmb{\P}}}(\pmb{\E},q_{\BS{\kx}_2}^*(\ko(a_1-2,b_1-2)))} \
commutes with base change. It is a vector bundle on $\pmb{\P}$ of rank
\m{h^1(\ko(a_0-a_1,b_0-b_1))}.

Let
\[\pmb{\M} \ = \ \P\big(\Ext^1_{q_{\pmb{\P}}}(\pmb{\E},
q_{\BS{\kx}_2}^*(\ko(a_1-2,b_1-2)))\big) \ , \]
which parameterizes non trivial extensions on \m{\BS{\kx}_2}
\[0\lra\ko(a_1-2,b_1-2)\lra\ke\lra\E\lra 0 \ , \]
with \m{\E} such that there is a non trivial extension $(\ref{equ10})$. We have
\[\dim(\pmb{\M}) \ = \ 2h^0(\ko(a_0-a_1,b_0-b_1))-2 \ = \
-2(a_1-a_0-1)(b_1-b_0-1) \ . \]
Let \
\m{\pmb{\M}'\subset\pmb{\M}} \ be the open subset corresponding to extensions
such that the associated morphism \ \m{\E\to\ko(a_1,b_1)} \ is surjective.

Let \ \m{\pi:\pmb{\M}\to\pmb{\P}}, \m{r_{\BS{\kx}_2}:
\BS{\kx}_2\times\pmb{\M}\to\BS{\kx}_2} \ and \
\m{r_{\pmb{\M}}:\BS{\kx}_2\times\pmb{\M}\to\pmb{\M}} \ be the projections.
There exists a {\em universal extension} on \ \m{\BS{\kx}_2\times\pmb{\M}}
\[0\lra r_{\BS{\kx}_2}^*(\ko(a_1-2,b_1-2))\ot
r_{\pmb{\M}}^*(\ko_{\pmb{\M}}(1))\lra\kh\lra\pi^\sharp(\E)\lra 0 \ . \]
Then \m{\pmb{\M}'}, equipped with \m{\kh_{|\BS{\kx}_2\times\pmb{\M}'}}, is a
fine moduli space for families of sheaves in \m{\Sigma(a_0,b_0,a_1,b_1)}
parameterized by integral schemes.

For every \m{m\in\pmb{\M}'}, the Koda\"ira-Spencer map \
\m{T_{\pmb{\M},m}\to\Ext^1_{\BS{\kx}_2}(\kh_m,\kh_m)} \ is an isomorphism \
\m{T_{\pmb{\M},m}\simeq H^1(\BS{\kx}_2,\EEnd(\kh_m))}.

\end{subsub}

\end{sub}

\vskip 4cm

\end{document}